\documentclass{article}
\usepackage{amsmath}
\usepackage{amsthm}
\usepackage{amssymb}
\usepackage{amsfonts}
\usepackage{mathrsfs}
\usepackage{enumerate}
\usepackage{authblk}
\usepackage[normalem]{ulem} 
\usepackage[dvipsnames]{xcolor}
\usepackage[utf8]{inputenc}
\usepackage[english]{babel}
\usepackage[ top = 3 cm, bottom = 2 cm]{geometry} 
\usepackage{graphicx}

\usepackage[colorlinks,citecolor=red,pagebackref,hypertexnames=false]{hyperref}
\usepackage{pdfsync}

\newcommand{\ep}{\varepsilon}
\newcommand{\n}[1]{\mathscr{#1}}
\newcommand{\m}[1]{\mathcal{#1}}
\newcommand{\bb}[1]{\mathbb{#1}}
\newcommand{\ra}{\rightarrow}

\newcommand{\D}{\Delta}
\DeclareMathOperator*{\esssup}{esssup}

\numberwithin{equation}{section}
\theoremstyle{plain}
\newtheorem{theorem}{Theorem}
\newtheorem{lemma}[theorem]{Lemma}
\newtheorem{proposition}[theorem]{Proposition}
\newtheorem{corollary}[theorem]{Corollary}

\theoremstyle{definition}
\newtheorem{definition}[theorem]{Definition}
\newtheorem{remark}[theorem]{Remark}

\begin{document}
\title{Optimal Stefan Problem}
\author[1]{Ugur G. Abdulla\thanks{abdulla@fit.edu}}
\author[2]{Bruno Poggi\thanks{poggi008@umn.edu}}
\affil[1]{Department of Mathematical Sciences, Florida Institute of Technology, Melbourne, Florida 32901}
\affil[2]{Department of Mathematics, University of Minnesota, Minneapolis, Minnesota 55455}
\maketitle
\abstract{We consider the inverse multiphase Stefan problem with homogeneous Dirichlet boundary condition on a bounded Lipschitz domain, where the density of the heat source is unknown in addition to the temperature and the phase transition boundaries. The variational formulation is pursued in the optimal control framework, where the density of the heat source is a control parameter, and the criteria for optimality is the minimization of the $L_2-$norm declination of the trace of the solution to the Stefan problem from a temperature measurement on the whole domain at the final time. The state vector solves the multiphase Stefan problem in a weak formulation, which is equivalent to Dirichlet problem for the quasilinear parabolic PDE with discontinuous coefficient. The optimal control problem is fully discretized using the method of finite differences. We prove the existence of the optimal control and the convergence of  the discrete optimal control problems to the original problem both with respect to cost functional and control. In particular, the convergence of the method of finite differences for the weak solution of the multidimensional multiphase Stefan problem is proved. The proofs are based on achieving a uniform $L_{\infty}$ bound and $W_2^{1,1}$ energy estimate for the discrete multiphase Stefan problem.}\\



\newpage
\section{Introduction}\label{description of results}
\subsection{Introduction and Motivation}\label{E:1:1}
Let $d\in\bb N, \Omega\subset \bb R^d$ be a bounded domain with Lipschitz boundary, $T>0$, and $D:=\Omega\times(0,T]$. Consider the general multi-dimensional multi-phase Stefan problem \cite{LSU}: given phase transition temperatures $u^1<u^2<\cdots<u^J$, find a temperature function $u:D\ra\bb R$ and the phase transition boundaries
\begin{equation}\nonumber
S^j=\{(x,t)\in D~|~u(x,t)=u^j\},\quad j=1,2,\ldots, J
\end{equation}
which satisfy
\begin{equation}\label{pde}
\alpha(u)u_t-\text{div}(k(u)\nabla u)=f(x,t),\qquad (x,t)\in D, \quad u(x,t)\neq u^j,  j=\overline{1,J}
\end{equation}
where $f$ is a known function, $\alpha,k$ are known positive functions which are smooth on each of the intervals $[u^j,u^{j+1}]$ and have discontinuities of the first kind at the points $u=u^j, j=1,\ldots, J$; 
\begin{align}
&[u]\big|_{S^j}=0, &\qquad j=\overline{1,J}, \label{ucont} \\[3mm]
&b_j\cos(\mathbf n,t)+\sum\limits_{i=1}^d[k(u)u_{x_i}]\cos(\mathbf n,x_i)\big|_{S^j}=0, &j=\overline{1,J}, \label{stefancond} \\[3mm]
&u(x,0) = \phi(x), &x\in\Omega\label{initial}, \\[3mm]
&u|_{S}=0,\label{flux}
\end{align}
where $\phi$ is a known function, each $b_j$ is a positive number, $\mathbf n$ is the normal to the free boundary $S^j$ in the direction of increasing $u$ (that is, along the gradient of $u$), and the saltus $[v]\big|_{S^j}$ is the difference between the limiting value of $v$ on $S^j$ when approached from the domains $\{(x,t)~|~u<u^j\}$ and $\{(x,t)~|~u>u^j\}$ respectively; $S=\partial\Omega\times(0,T]$ is a lateral boundary of the cylinder $D$.

In the physical context, $f$ characterizes the density of the sources, $\phi$ is the initial temperature, (\ref{stefancond}) is the Stefan condition expressing the conservation law according to which the free boundary is pushed by the saltus of the heat flux from different phases, and (\ref{flux}) states that the temperature at the boundary is held constant at $0$.

Weak formulation of the multiphase Stefan problem, as well as existence and uniqueness of the weak solution to the multiphase Stefan problem was first proved in \cite{Kamenomostskaya, Oleinik}. We refer to monographies \cite{LSU, Meyrmanov} for the extensive list of references. 

Assume now that some of the data is not available, or involves some measurement error. For example, suppose that the density of the heat sources $f$ is not known and must be found along with the temperature $u$ and the free boundaries $S^j$. As compensation for not knowing this function, we must have access to additional information, which for instance may come as a measurement of the temperature at the final moment:
\begin{equation}\label{extra}
u\big|_{\Omega\times\{t=T\}}=\nu.
\end{equation}

\textbf{Inverse Multiphase Stefan Problem (IMSP).} Find the temperature function $u(x,t)$, free boundaries $S^j, j=1,...,J$, and the density of the heat sources $f(x,t)$ satisfying (\ref{pde})-(\ref{extra}). \\

The IMSP is not well posed in the sense of Hadamard. That is, if the data is not sufficiently coordinated, there may be no solution. Even if it exists, it might be not unique, and most importantly there is in general no continuous dependence of the solution on the data functions.  

In two recent papers \cite{Abdulla1, Abdulla2} a new variational formulation of the one-phase inverse Stefan problem (ISP) was developed when space dimension is one. An optimal control framework was implemented in which the boundary heat flux and the free boundary are components of the control vector and the optimality criteri consists of the minimization of the sum of $L_2$-norm declinations from the available measurement of the temperature on the fixed boundary and available 
information on the phase transition temperature on the free boundary. This approach allows 
one to tackle situations when the phase transition temperature is not known explicitly, and is available through measurement with possible error. It also allows for the development of iterative numerical methods of least computational cost due to the fact that for every given control vector, the parabolic PDE is solved in a fixed region instead of full free boundary problem. In \cite{Abdulla1} the well-posedness in Sobolev spaces framework and 
convergence of time-discretized optimal control problems is proved. In \cite{Abdulla2} full discretization was implemented and the convergence of the discrete optimal control problems to the original problem both with respect to cost functional and control is proved.  The main advantage of this method is that numerically, the problem to be solved at each step is only a Neumann problem, and not a full free boundary problem. In \cite{Abdulla3,Abdulla4} Frechet differentiability and first order optimality condition in Besov spaces framework is proved and the formula for the Frechet gradient is derived. Numerical analysis via iterative gradient method in Hilbert-Besov spaces based on the results of \cite{Abdulla1,Abdulla2,Abdulla3,Abdulla4} was implemented in \cite{Abdulla5}.

The new variational approach developed in \cite{Abdulla1,Abdulla2} is not applicable to the multiphase Stefan problem. The reason is that the Stefan condition on the phase transition boundary includes the flux calculated from both phases. Therefore, it can't be treated as a Neumann condition, even if we include the free boundary as one of the control components. In \cite{PoggiAbdulla} a new approach was developed based on the weak formulation of the multiphase Stefan problem as a boundary value problem for the nonlinear PDE with discontinuous coefficient. The optimal control framework was applied to the inverse multiphase Stefan problem with non-homogeneous Neumann conditions on the fixed boundaries in the case when the space dimension is one. The control vector was taken to be the heat flux on the left boundary and the optimality criteria consisted of the $L_2-$norm declinations from a measurement of the temperature on the right fixed boundary. The full discretization was implemented and convergence of the discrete optimal control problems to the original problem was proved. 

The main goal of this paper is to apply the idea of the paper \cite{PoggiAbdulla} to IMSP when the number of spatial dimensions is larger than $1$. We prove the existence of the optimal control and convergence of the sequence of discrete optimal control problems to the continuous problem both with respect to the functional and control. The proof is based on the proof of uniform $L_{\infty}$ bound, and $W_2^{1,1}$-energy estimate for the discrete multiphase Stefan problem, and results on the convergence of suitable interpolations of the discrete solutions. We address the problem of Frechet differentiability and application of the iterative gradient methods in Hilbert spaces in an upcoming paper.

We refer to a recent paper \cite{Abdulla1} for review of the literature on Inverse Stefan Problems. Most of the papers on ISP are in the one-dimensional case. Inverse Stefan problems with given phase boundaries were considered in \cite{Alifanov,Bell,Budak,BudakVasileva1, BudakVasileva2,Cannon,Cannon3,Carasso,Ewing1,Ewing2,Hoffman,Sherman,Goldman}; optimal control of Stefan problems, or equivalently inverse problems with unknown phase boundaries were investigated in \cite{Baumeister,Fasano,Hoffman1,Hoffman2,Jochum2,Jochum1,Knabner,Lurye,Nochetto, Niezgodka,Primicero,Sagues,Talenti,Goldman,Gol'dman,Vasilev,Yurii}. 

The structure of the paper is as follows: in Section~\ref{E:1:1a} the notation of Sobolev spaces are described. In Section~\ref{E:1:1b}  we formulate the IMSP as an optimal control problem. In Section~\ref{E:1:1c} we perform full discretization through finite differences and formulate discrete optimal control problem. In Section \ref{assumptions}, all the operative assumptions are declared. In Section~\ref{mainresults} the main results are formulated.  In Section~\ref{prelim} we prove the existence and uniqueness of the discrete state vector, as well as other auxiliary lemmas. In Section~\ref{estimates}, we prove $L_{\infty}$ and $W_2^{1,1}$ estimates that the discrete state vectors satisfy. Section~\ref{Interpolations} describes different interpolations of the discrete state vectors to the whole domain and contains proofs on appropriate equivalences of the different interpolations. In Section~\ref{approxtheorem}, it is shown that piece-wise linear interpolations approximate a weak solution to the Stefan problem. This allows us to prove in Section~\ref{exists} the existence of a solution to the optimal control problem, and in Section~\ref{approxfunc} we prove convergence of the discrete optimal control problems to the continuous optimal control problem.

\subsection{Notations}\label{E:1:1a}
$B_r(x)\subset \bb R^d$ - ball of radius $r$ and center $x$; $m_d(\cdot)$ - $d$-dimensional Lebesgue measure; 
\[
\Omega+z:=\{x\in\bb R^d~|~\exists y\in\Omega\text{ s.t. }y+z=x\},\qquad\Omega+A:=\bigcup\limits_{z\in A}(\Omega+z) \ \text{for} \ A\subset\bb R^d
\]
$L_p(D), 1\leq p <+\infty$ - Banach space of real-valued measurable functions on $D$ with finite norm
\[
\Vert f\Vert_{L_p(D)}:=\Big (\int\limits_{D}|f|^p\,dx\Big )^{\frac{1}{p}}<+\infty,
\]
$L_{\infty}(D)$ - Banach space of essentially bounded real-valued measurable functions on $D$ with norm
\[
\Vert f\Vert_{L_{\infty}(D)}=\esssup\limits_{(x,t)\in D}|f(x,t)|<+\infty.
\]
$W_2^1(\Omega)$ -  Hilbert space of all elements $f$ of $L_2(\Omega)$ for which the partial weak derivative $\partial f/\partial x_i$ exists and lie in $L_2(\Omega)$ for each $i=1,\ldots,d$. This space has inner product
\[
(f,g)=\int\limits_{\Omega}\left(fg+\sum\limits_{i=1}^d\frac{\partial f}{\partial x_i}\frac{\partial g}{\partial x_i}\right)\,dx.
\]
$W_2^{1,0}(D)$ -  Hilbert space of all elements $f$ of $L_2(D)$ having square-integrable first-order weak partial derivatives in all spatial directions. This space is endowed with the inner product
\[
(f,g)=\int\limits_{D}\left(fg+\sum\limits_{i=1}^d\frac{\partial f}{\partial x_i}\frac{\partial g}{\partial x_i}\right)\,dx\,dt.
\]
$W_2^{1,1}(D)$ - Hilbert space of all elements of $L_2(D)$ having square-integrable first-order weak partial derivatives in all coordinate directions. The inner product is
\[
(f,g)=\int\limits_{D}\left(fg+\frac{\partial f}{\partial t}\frac{\partial g}{\partial t}+\sum\limits_{i=1}^d\frac{\partial f}{\partial x_i}\frac{\partial g}{\partial x_i}\right)\,dx\,dt.
\]
$\overset{\circ}W{}_2^{1,1}(D)$ - linear subspace of elements $f$ of $W_2^{1,1}(D)$ which satisfy
\[
f\Big|_S=0,
\]
in the sense of traces.

\subsection{Multiphase Stefan Optimal Control Problem}\label{E:1:1b}

Following the usual reformulation of the inverse multiphase Stefan problem (see \cite{LSU, Oleinik}), we define the function
\begin{equation}\label{F}
F(t)=\int\limits_{0}^tk(y)\,dy,
\end{equation}
and consider the transformation
\begin{equation}\label{vt}
v(x,t):=F(u(x,t)).
\end{equation}
Then $v^j=F(u^j)$, $~v^1<\cdots<v^J$, and our conditions become:

\begin{align}
&\beta(v)v_t-\Delta v=f(x,t), &\qquad (x,t)\in D, v(x,t)\neq v^j,\label{bvpde}\\
&[v]|_{S^j}=0, &\qquad j=\overline{1,J},\label{vcont} \\
&b_j\cos(\mathbf n,t)+\sum\limits_{i=1}^d[v_{x_i}]\cos(\mathbf n,x_i)\big|_{S^j}=0, &\qquad j=\overline{1,J},\label{vjump} \\
&v\big|_{\Omega\times\{t=0\}}=\Phi:=F(\phi),&\label{vphi}\\
&v|_{S}=0,&  \label{vg}\\
&v\big|_{\Omega\times\{t=T\}}=\Gamma:=F(\nu),&\label{vnu}
\end{align}
with $\beta(v)$ possessing the same properties as $\alpha(u)$. Now, we can invoke a monotone increasing piecewise smooth function $b(v)$ such that $b'(v)=\beta(v)$ on each of the intervals $(v^j,v^{j+1})$.  Our partial differential equation becomes
\begin{equation}\label{bpde}
\frac{\partial b(v)}{\partial t}-\Delta v=f(x,t), \qquad (x,t)\in D,~ v(x,t)\neq v^j.
\end{equation}

Moreover, we're free to choose the jump of $b$ at the values $v=v^j$. We choose them in such a way that $[b(v)]|_{S^j}=-b_j$ so that upon integration by parts of (\ref{bpde}) over $D$, the integrals over the phase transition boundaries cancel out.

\begin{definition}\label{typeB} We say that a measurable function $B(x,t,v)$ is \textit{of type }$\n B$ if
\begin{enumerate}[(a)]
	\item $B(x,t,v)=b(v),\qquad v\neq v^j,\quad\forall j=\overline{1,J}$
	\item $B(x,t,v)\in[b(v^j)^-,b(v^j)^+],\qquad v=v^j$ for some $j$. \\
\end{enumerate}
Note that $B(x,t,v)$ can take different values for different $(x,t)$ when $v=v^j$ for some $j$.
\end{definition}

Given $f$, a solution to the Stefan problem (\ref{bvpde})-(\ref{vg}) is understood in the following sense:

\begin{definition}\label{weaksoldef} $v\in\overset{\circ}{W}{}_2^{1,1}(D)\cap L_{\infty}(D)$ is called a \emph{weak solution of the Stefan problem}  (\ref{bvpde})-(\ref{vg}) if for any two functions $B,B_0$ of type $\n B$, the integral identity
\begin{gather}
\int\limits_D\Big[-B(x,t,v(x,t))\psi_t+\nabla v\cdot\nabla\psi-f\psi\Big]\,dxdt - \int\limits_{\Omega}B_0(x,0,\Phi(x))\psi(x,0)\,dx=0\label{weaksol}
\end{gather}
is satisfied for arbitrary $\psi\in\overset{\circ}{W}{}_2^{1,1}(D)$ with { $\psi|_{\Omega\times\{t=T\}}=0$}.
\end{definition}

For fixed $R>0$, define the \emph{continuous control set}
\[\n F^R =\left\{f\in L_{\infty}(D)~\Bigg|~\Vert f\Vert_{L_{\infty}(D)}\leq R\right\}.\]
Consider minimization of the cost functional 
\begin{equation}\label{functional}
\n J(f)=\Vert v|_{\Omega\times\{t=T\}}-\Gamma\Vert^2_{L_2(\Omega)}
\end{equation}
on $\n F^R$, where $v=v(x,t;f)\in\overset{\circ}{W}{}_2^{1,1}(D)\cap L_{\infty}(D)$ is a weak solution of the Stefan problem in the sense of Definition \ref{weaksoldef}. This optimal control problem will be called \emph{Problem $\m I$}.

\subsection{Discrete Optimal Control Problem}\label{E:1:1c}

We apply the method of finite differences. 
Let $n\in\bb N, \tau:=\frac Tn, h>0$, and cut $\bb R^d\times\bb R$ by the planes 
\[
x_i=k_ih, ~i=1,\ldots,d,\quad t=k_0\tau,\qquad\forall k_{\ell}\in\bb Z,~\ell=0,1,\ldots,d,
\]
so as to obtain a collection of elementary (closed) cells with length $h$ in each $x_i$ direction and length $\tau$ in the $t$ direction. { We will denote by $\Delta$ the discretization with steps $(\tau,h)$. We introduce a partial ordering on the set of discretizations: we say that $\D_1\leq\D_2$ if $\tau_1\leq\tau_2$ and $h_1\leq h_2$. We will call $t_{\ell}=\tau\ell$ for $\ell=1,\ldots, n$.} Let $\alpha=(k_1,k_2,\ldots,k_d,k_0)$ be a multi-index, and $\gamma=(k_1,k_2,\ldots,k_d)$. { We will agree to write $\alpha=(\gamma,k_0)$, $\alpha_i$ is the $i-$th component of $\alpha$ if $i\in\{1,2,\ldots, d\}$ and $\alpha_0$ is the $d+1-$st component of $\alpha$, while $\gamma_i$ is the $i-$th component of $\gamma$}. Then each elementary cell $C^{\alpha}_{\Delta}$ can be written uniquely in the following way
\[
C^{\alpha}_{\Delta}=\Big\{(x,t)\in\bb R^d\times\bb R~\big|~k_ih\leq x_i\leq(k_i+1)h,~i=1,\ldots,d;~~(k_0-1)\tau\leq t\leq k_0\tau\Big\}.
\]

Similarly we define the rectangular prisms:
\[
R_{\D}^{\gamma}=\Big\{x\in\bb R^d~\big|~k_ih\leq x_i\leq(k_i+1)h,~i=1,\ldots,d\}.
\]
and whenever we write $k$ as a superscript to a set in $\bb R^d$, it is meant the projection of that set onto the hyper-plane $t=k\tau$ of $R^{d+1}$. For instance,
\[
R^{\gamma,k}_{\D}=\Big\{(x,t)\in\bb R^d\times\bb R~\big|~k_ih\leq x_i\leq(k_i+1)h,~i=1,\ldots,d;~~t=k\tau\}.
\]
We write the collections of these cells and prisms as
\[
\n C_{\D}=\big\{ C^{\alpha}_{\D}~|~\alpha\in\bb Z^{d+1}\big\},
\]
\[
\n R_{\D}=\big\{ R^{\gamma}_{\D}~|~\gamma\in\bb Z^{d}\big\},
\]
and consider the subcollections which lie only in $\overline D$ and $\overline\Omega$ respectively:
\[
\n C^D_{\D}=\Big\{C^{\alpha}_{\D}\in\n C_{\D}~|~C_{\D}^{\alpha}\subset\overline D\Big\},
\]
\[
\n R^{\Omega}_{\D}=\Big\{R^{\gamma}_{\D}\in\n R_{\D}~|~R_{\D}^{\gamma}\subset\overline{\Omega}\Big\}.
\]
The unions of the elements in these subcollections comprise the discretized versions of $D$ and $\Omega$ respectively. So we write
\[
\Omega_{\D}=\bigcup_{R^{\gamma}_{\D}\in\n R^{\Omega}_{\D}}R^{\gamma}_{\D}~\subset\overline{\Omega},\qquad D_{\D}=\bigcup_{C^{\alpha}_{\D}\in\n C^D_{\D}}C^{\alpha}_{\D}~\subset\overline D.
\]
By the \emph{natural corner} of a prism in $\n R_{\D}$ it is meant the vertex of the prism whose coordinates are smallest relative to the other vertexes, and by the natural corner of a cell $C_{\D}^{(\gamma,k)}\in\n C_{\D}$ it is meant the vertex of the cell whose spatial coordinates are the same as those of the natural corner of $R_{\D}^{\gamma}$, and whose time coordinate is $k\tau$. From here on, we identify each prism (cell) by its natural corner.\\

We denote by $S_{\D}$ the lateral boundary of $D_{\D}$, $D_{\D}'= (D_{\D}\backslash\partial D_{\D})\cup(\Omega_{\D}\times\{t=T\})$ and $\Omega_{\D}^{'}=\Omega_{\D}\backslash\partial\Omega_{\D}$. Now define the lattice of points

\[
\n L_T=\Big\{(x,t)\in\bb R^d\times\bb R~|~\exists\alpha\in\bb Z^{d+1}\text{ s.t. }x_i=k_ih,~i=1,\ldots,d,~~t=k_0\tau\Big\},
\]
\[
\n L=\Big\{x\in\bb R^d~|~\exists\gamma\in\bb Z^{d}\text{ s.t. }x_i=k_ih,~i=1,\ldots,d\Big\}.
\]
We will usually write {$y=(x,t)$, $y_{\alpha}=(k_1h,k_2h,\ldots,k_dh,k_0\tau), x_{\gamma}=(k_1h,k_2h,\ldots,k_dh)$}. Note the obvious bijections $\alpha\mapsto y_{\alpha}$, $\gamma\mapsto x_{\gamma}$; bijections of this form will henceforth be referred as natural. Given a set $X$ which is in natural bijection with a subset of the set of multi-indexes $\gamma$ (or $\alpha$), we write $\n A(X)$ as the indexing set. Moreover, if $X\subset\bb R^d$, then $\n L(X):=\n L\cap X$ (and similarly if $X\subset\bb R^{d+1}$). When $X=\n L(Y)\subset\bb R^d$, we'll agree to write $\n A(Y)$ instead of $\n A(\n L(Y))$ (and likewise if $X=\n L_T(Y)$). For emphasis, by $\n A:=\n A(\n R^{\Omega}_{\D})$ it is meant the set of all those indexes $\gamma$ which correspond to a prism in $\Omega_{\D}$. These indexes are also in natural bijection with the natural corners of these prisms. In particular, some of the corresponding lattice points may fall on the boundary $\partial\Omega_{\D}$. We contrast this set to the set $\n A(\Omega_{\D}')$ of indexes in natural bijection to the lattice points that lie strictly in the interior of $\Omega_{\D}$, and to the set $\n A(\Omega_{\D})$, of all indexes which are in natural bijection with the lattice points that lie in $\Omega_{\D}$. It is clear that $\n A(\Omega_{\D}')$ is a subset of $\n A$. For ease of notation, we will often write
\[
\sum\limits_{\n A(X)}\quad\text{instead of}\quad\sum\limits_{\gamma\in\n A(X)},
\]
and likewise for other expressions requiring subscripts.\\

It will be important to give a sense as to how to discretize functions given in the continuous setting.  Given $\Phi\in W_2^1(\Omega),\Gamma\in L_2(\Omega), f\in L_2(D)$, we will construct appropriately discretized versions of these functions through the use of the Steklov averages. { First fix an extension of $\Phi$ to $\Omega+B_1(0)$ so that the extension lies in $W_2^1(\Omega+B_1(0))$. Henceforth refer to the extension as $\Phi$. We denote \begin{equation}\label{phid}
\Phi_{\gamma}=\frac1{h^d}\int\limits_{x_1}^{x_1+h}~\int\limits_{x_2}^{x_2+h}\cdots\int\limits_{x_d}^{x_d+h}\Phi(x)\,dx,\qquad\text{where }\gamma\in\n A(\Omega_{\D}),~~ x_i\text{ is }i\text{-th coordinate of }x_{\gamma},
\end{equation}
and
\begin{equation}\nonumber
\Gamma_{\gamma}=\frac1{h^d}\int\limits_{x_1}^{x_1+h}~\int\limits_{x_2}^{x_2+h}\cdots\int\limits_{x_d}^{x_d+h}\Gamma(x)\,dx,\qquad\text{where }\gamma\in\n A,~~ x_i\text{ is }i\text{-th coordinate of }x_{\gamma}.
\end{equation}}
We note the region of integration in (\ref{phid}) is $R_{\D}^{\gamma}$. Also,
\begin{equation}\label{fd}
f_{\alpha}=\frac1{\tau h^d}\int\limits_{t_{k-1}}^{t_k}~\int\limits_{x_1}^{x_1+h}~\int\limits_{x_2}^{x_2+h}\cdots\int\limits_{x_d}^{x_d+h}f(x,t)\,dx\,dt,\qquad\alpha=(\gamma,k)\in\n A(\n C_{\D}^D),
\end{equation}
and we observe the region of integration in (\ref{fd}) is really $C_{\D}^{\alpha}$.

We will need to smoothen the function $b$. To this end, for $\rho>0$ let
$\omega_{\rho}$ be a non-negative $C_0^{\infty}(\bb R)$ mollifier. We can take, for example,
\begin{equation}\label{kernel}
\omega_{\rho}(v) =\left\{\begin{matrix}\m C \rho^{-1}e^{-\frac{\rho^2} {\rho^2-v^2}},\quad&|v|\leq\rho\\ 0,\quad&|v|>\rho\end{matrix}\right.
\end{equation}
where $\m C$ is a constant chosen so that $\int\limits_{\bb R}\omega_1(|u|)\,du =1$. We then define
\begin{equation}\label{bn}
b_n:=b*\omega_{\frac1n}.
\end{equation}

Given a discretization $\D$, we use the notation $[f]_{\D}$ for a collection of real numbers $\{f_{\alpha}\},~\alpha\in\n A(\n C_{\D}^D)$. Each of these can be thought of as vectors in a suitable finite-dimensional space. We define
\begin{equation}\nonumber
\Vert[f]_{\D}\Vert_{\ell_{\infty}}:=\max\limits_{\n A(\n C_{\D}^D)}|f_{\alpha}|, \ \   \Vert[f]_{\D}\Vert_{\ell_2}:=\Big (\sum\limits_{\n A(\n C_{\D}^D)}\tau h^d f_{\alpha}^2\Big )^{\frac{1}{2}}.
\end{equation}
We will consider space and time differences. For a collection of numbers $\{u_{\alpha}\}$, if we write $\alpha=(\gamma, k_0)$, then
\[
u_{\alpha\bar t}=\frac{u_{(\gamma,k_0)}-u_{(\gamma,k_0-1)}}{\tau}.
\]
is the backward time difference. The forward space difference along the $x_i-$direction $u_{\alpha x_i}$ is
\[
u_{\alpha x_i}=\frac{u_{(k_1,\ldots,k_i+1,\ldots,k_d,k_0)}-u_{(k_1,\ldots,k_i,\ldots,k_d,k_0)}}{h}.
\]
Moreover, for convenience of notation, we will write
\[
\gamma+e_i:=(k_1,\ldots,k_i+1,\ldots,k_d),\quad \alpha+e_i:=(k_1,\ldots,k_i+1,\ldots,k_d,k_0)
\]
for suitable $i$. For fixed $R>0$, define the \emph{discrete control sets}
\[
\n F_{\D}^R:=\Big\{[f]_{\D}~\big|~\Vert[f]_{\D}\Vert_{\ell_{\infty}}\leq R\Big\}
\]
and the following mappings between the continuous and discrete control sets. Let \[
\n P_{\D}:\bigcup_R\n F_{\D}^R\longrightarrow\bigcup_R\n F^R,\qquad\n P_{\D}([f]_{\D})=f^{\D}
\]
be an interpolating map, where
\[
f^{\D}\Big|_{C_{\D}^{\alpha}}=f_{\alpha},~~\alpha\in\n A(\n C_{\D}^D),\qquad f^{\D}\equiv0~~\text{elsewhere on }D.
\]
Also, let
\[
\n Q_{\D}:\bigcup_R\n F^R\longrightarrow\bigcup_R\n F_{\D}^R,\qquad\n Q_{\D}(f)=[f]_{\D}
\]
be a discretizing map, where $f_{\alpha}$ is given by (\ref{fd}) for each $\alpha\in\n A(\n C_{\D}^D)$.\\

At this point we are ready to define a solution of the discrete Stefan problem. 

\begin{definition}\label{dsvdef} Given $[f]_{\D}$, the vector function $[v([f]_{\D})]_{\D}=(v(0),v(1),\ldots,v(n))$, $v(k)$ a collection of real numbers $\{v_{\gamma}(k)\}, \gamma\in\n A(\Omega_{\D}),~k=0,1,2\ldots,n$, is called a \emph{discrete state vector} provided it satisfies
\begin{enumerate}[{\bf(i)}]
\item $v_{\gamma}(0)=\Phi_{\gamma},~\gamma\in\n A(\Omega_{\D}')$,
\item For each fixed $k=1,\ldots,n$, the collection $v(k)$ satisfies
\begin{align}
\sum\limits_{{\n A}}h^d\Bigg[\Big(b_n(v_{\gamma}(k))\Big)_{\bar t}\eta_{\gamma}~+ ~&\sum\limits_{i=1}^dv_{\gamma x_i}(k)\eta_{\gamma x_i}~-f^{\D}_{(\gamma,k)}\eta_{\gamma}\Bigg]=0\label{dsveq}
\end{align}
for arbitrary collection of values $\{\eta_{\gamma}\},~\gamma\in\n A(\Omega_{\D})$ which satisfies that $\eta_{\gamma}=0$ for $\gamma\in\n A(\partial\Omega_{\D})$.
\item For each $k=0,1,\ldots,n$, we have $v_{\gamma}(k)=0$ for $\gamma\in\n A(\partial\Omega_{\D})$.
\end{enumerate}
\end{definition}

We note that the collection $\{f^{\D}_{\alpha}\}$ appearing in (\ref{dsveq}) is the function $\n Q_{\D}\Big(\n P_{\D}([f]_{\alpha})\Big)$. Given $[f]_{\D}\in\n F_{\D}^R$ for some $R>0$, it will be shown that the discrete state vector $[v([f]_{\D})]_{\D}$ exists uniquely. This allows us to define a discrete cost functional $\n I_{\D}:\cup_R\n F_{\D}^R\ra[0,+\infty)$ by
\begin{equation}\label{discretecostfunctional}
\n I_{\D}([f]_{\D})=\sum\limits_{\n A}h^d|v_{\gamma}(n)-\Gamma_{\gamma}|^2
\end{equation}
where the $v_{\gamma}(n)$ are taken from $v(n)$, the $n-$th component of the discrete state vector $[v([f]_{\D})]_{\D}$. The discrete optimal control problems will be called problems $\m I_{\D}$. We define
\begin{equation}\label{zeta1}
\zeta_{\D,k}^{\gamma}:=\int\limits_0^1b_n'\Big(\theta v_{\gamma}(k)+(1-\theta)v_{\gamma}(k-1)\Big)\,d\theta,
\end{equation}
for each $(\gamma,k)\in\n A(D_{\D}), k\neq0$. For each such $(\gamma,k)$ we note
\begin{equation}\label{mvt}
b_n(v_{\gamma}(k))-b_n(v_{\gamma}(k-1))=\zeta_{\D,k}^{\gamma}\Big(v_{\gamma}(k)-v_{\gamma}(k-1)\Big).
\end{equation}

\subsection{Assumptions}\label{assumptions}

Throughout the paper we will make the following assumptions:
\begin{enumerate}[(a)]
	\item $\Omega\subset\bb R^d$ is open, bounded, and has Lipschitz boundary.
	\item $\alpha$ and $k$ are positive on $\bb R$, and the restrictions of $\alpha,k$ to each of the segments $(-\infty,u^1),$ $ (u^j,u^{j+1}), j=1,\ldots,J-1,~(u^J,+\infty)$ are continuously differentiable functions with positive limits at the finite end-points.
	\item $\min\left\{\liminf\limits_{u\ra+\infty}\frac{\alpha(u)}{k(u)}~,~\liminf\limits_{u\ra-\infty}\frac{\alpha(u)}{k(u)}\right\}\geq a_0$ for some $a_0\in(0,+\infty)$.
	\item
	\begin{equation}\label{infiniteint}
	\int\limits_{0}^{+\infty}k(y)\,dy=+\infty,\qquad\int\limits_{0}^{-\infty}k(y)\,dy=-\infty.
	\end{equation}
	\item $\phi\in W_2^1(\Omega)\cap L_{\infty}(\Omega)$.
	\item Either $\nu\in L_{\infty}(\Omega)$, or $\nu\in L_2(\Omega)$ and $k\in L_{\infty}(\bb R)$.
	\item For each $j=1,\ldots,J$, the set $\big\{x~|~\phi(x)=u^j\big\}$ has $d-$dimensional measure $0$.
\end{enumerate}

A brief discussion of the assumptions follows: Assumption that $\Omega$ is Lipschitz is assumed to guarantee application of standard Sobolev embedding theorems. Assumption (b) allows for the function $b$ to be continuously differentiable and strictly monotone increasing on each of the segments $(-\infty,u^1),(u^j,u^{j+1}),j=1,\ldots,J-1,(u^J,+\infty)$. Assumption (c) provides positive lower bound for $b$, which we use to prove the existence of the discrete state vector, as well as to establish the energy estimates. Given assumption (b), assumption (d) is a necessary and sufficient condition that the map $F:\bb R\ra\bb R$ is a bijection. 
Furthermore, assumption (d) allows the function $b$ to have the aforementioned properties on all of $\bb R$, a requisite for our proof of the existence of the discrete state vector. Assumption (e) is important for the energy estimates. Either of the conditions in assumption (f) will guarantee $\Gamma\in L_2(\Omega)$, which allows for the functional $\n J$ to be well-defined. Finally, assumption (g) guarantees that the second term in the integral identity \eqref{weaksol} is independent of the choice of the functions $B_0$ of type $\n B$.

\subsection{Main Results}\label{mainresults}

We have the following results:

\begin{theorem}\label{optsol} The optimal control problem $\m I$ has a solution. That is, the set
	\[
	\n F_*:=\Big\{f\in\n F^R~\Big|~\n J(f)=\n J_*:=\inf\limits_{f\in\n F^R}\n J(f)\Big\}
	\]
	is not empty.
\end{theorem}

\begin{theorem}\label{funcapprox}The sequence of discrete optimal control problems $\m I_n$ approximates the optimal control problem $\m I$ with respect to the functional, that is,
\begin{equation}\label{approx1}
\lim\limits_{\D\ra0}\n I_{\D_*}=\n J_*
\end{equation}
where
\[
\n I_{\D_*}=\inf\limits_{\n F_{\D}^R}\n I([f]_{\D}).
\]
	
Furthermore, let $\{\ep_{\D}\}$ be a sequence of positive real numbers with $\lim\limits_{\D\ra0}\ep_{\D}=0$. If the sequence $[f]_{\D,\ep}\in\n F_{\D}^R$ is chosen so that
\begin{equation}\label{approxcond}
\n I_{\D_*}\leq\n I_{\D}([f]_{\D,\ep})\leq\n I_{\D_*}+\ep_{\D},
\end{equation}
then we have
\begin{equation}\label{approx2}
\lim\limits_{\D\ra0}\n J(\n P_{\D}([f]_{\D,\ep}))=\n J_*.
\end{equation}
Also, the sequence $\{\n P_{\D}([f]_{\D,\ep})\}$ is uniformly bounded in $L_2(D)$ and all of its $L_2(D)-$weak limit points lie in $\n F_*$. Moreover, if $f_*$ is such a weak limit point, then there is a subsequence $\D'$ such that the linear interpolations $V_{\D'}'$ of the discrete state vectors $[v([f]_{\D',\ep})]_{\D'}$ converge weakly in $W_2^{1,1}(D)$ to $v=v(x,t;f_*)$, a weak solution to the Stefan Problem in the sense of Definition \ref{weaksoldef}.
\end{theorem}

\section{Preliminary Results}\label{prelim}

\begin{proposition}\label{equiv} Fix a discretization $\D$ and control $[f]_{\D}$. For a vector function $[v([f]_{\D})]_{\D}$ as in Definition \ref{dsvdef}, consider the following condition:
	
\textbf{(ii)'}$~$ For each $k=1,2,\ldots,n$ and $\gamma$ such that $\gamma\in\n A(\Omega_{\D}')$, we have
\begin{equation}\label{dsvequiv}
\Big(b_n(v_{\gamma}(k))\Big)_{\bar t}-\Delta_hv_{\gamma}(k)=f^{\D}_{(\gamma,k)},
\end{equation}
where
\begin{equation}
\Delta_h v_{\gamma}(k):=\sum\limits_{i=1}^d\Big(v_{\gamma x_i}(k)\Big)_{\overline{x_i}}.\nonumber
\end{equation}
Then $[v([f]_{\D})]_{\D}$ is a discrete state vector if and only if it satisfies conditions (i), (ii)', and (iii).
\end{proposition}

\emph{Proof.} Suppose $[v([f]_{\D})]_{\D}$ satisfies (i),(ii)' and (iii). Fix $k\in\{1,2,\ldots,n\}$. Consider an arbitrary collection $\{\eta_{\gamma}\}$ of real numbers for $\gamma\in\n A(\Omega_{\D})$ which satisfies $\eta_{\gamma}=0$ for $\gamma\in\n A(\partial\Omega_{\D})$. For each $\gamma\in\n A(\Omega_{\D}')$, multiply (\ref{dsvequiv}) by $h^d\eta_{\gamma}$, and then perform a summation of all (\ref{dsvequiv}) over $\gamma\in\n A(\Omega_{\D}')$. We obtain
\begin{equation}\label{almostdsveq}
\sum\limits_{\n A(\Omega'_{\D})}h^d\Bigg[\Big(b_n(v_{\gamma}(k))\Big)_{\bar t}\eta_{\gamma}~-\Delta_h v_{\gamma}(k)\eta_{\gamma}~-f^{\D}_{(\gamma,k)}\eta_{\gamma}\Bigg]=0.
\end{equation}

Observe that
\begin{gather*}
-\sum\limits_{\n A(\Omega'_{\D})}\sum\limits_{i=1}^d\Big(v_{\gamma x_i}(k)\Big)_{\overline{x_i}}\eta_{\gamma}=-\sum\limits_{\n A(\Omega'_{\D})}\sum\limits_{i=1}^d\frac{v_{\gamma x_i}(k)-v_{(\gamma-e_i)x_i}(k)}h\eta_{\gamma}=\\[4mm]=-\sum\limits_{\n A(\Omega'_{\D})}\sum\limits_{i=1}^d\frac{v_{\gamma x_i}(k)}h\eta_{\gamma}+\sum\limits_{\n A(\Omega'_{\D})}\sum\limits_{i=1}^d\frac{v_{(\gamma-e_i)x_i}(k)}h\eta_{\gamma},\\\intertext{letting $z=\gamma-e_i$ in the last summation we see}\\=-\sum\limits_{\n A(\Omega'_{\D})}\sum\limits_{i=1}^d\frac{v_{\gamma x_i}(k)}h\eta_{\gamma}+\sum\limits_{z\text{ and }i\text{ s.t. } z+e_i\in\n A(\Omega'_{\D})}\frac{v_{z x_i}(k)}h\eta_{z+e_i}=\\[4mm]=\sum\limits_{\gamma\in\n A(\Omega_{\D}')\text{ and }i\text{ s.t. } \gamma+e_i\in\n A(\Omega'_{\D})}v_{\gamma x_i}(k)\eta_{\gamma x_i}-\sum\limits_{\gamma\in\n A(\Omega_{\D}')\text{ and }i\text{ s.t. }\gamma+e_i\in\n A(\partial\Omega_{\D})}\frac{v_{\gamma x_i}(k)}h\eta_{\gamma}~+\\[4mm]+\sum\limits_{z\in\n A(\partial\Omega_{\D})\text{ and }i\text{ s.t. } z+e_i\in\n A(\Omega'_{\D})}\frac{v_{z x_i}(k)}h\eta_{z+e_i}=\\[4mm]=\sum\limits_{\gamma\in\n A(\Omega_{\D}'),i~|~ \gamma+e_i\in\n A(\Omega'_{\D})}v_{\gamma x_i}(k)\eta_{\gamma x_i}+\sum\limits_{\gamma\in\n A(\Omega_{\D}'),i~|~\gamma+e_i\in\n A(\partial\Omega_{\D})}\frac{v_{\gamma x_i}(k)}h\big(-\eta_{\gamma}+\eta_{\gamma+e_i}\big)~+\\[4mm]+\sum\limits_{z\in\n A(\partial\Omega_{\D})\text{ and }i\text{ s.t. } z+e_i\in\n A(\Omega'_{\D})}\frac{v_{z x_i}(k)}h\big(\eta_{z+e_i}-\eta_z\big)=\\[4mm] =\sum\limits_{\n A}\sum\limits_{i=1}^dv_{\gamma x_i}(k)\eta_{\gamma x_i}.
\end{gather*}
Plugging this calculation into (\ref{almostdsveq}) and using the fact that $\eta_{\gamma}=0$ for each $\gamma\in\n A\backslash\n A(\Omega_{\D}')$ shows that (ii) is satisfied. Conversely, suppose (i), (ii) and (iii) are satisfied, and fix $k\in\{1,2,\ldots,n\}$. Fix an arbitrary $\gamma'$ such that $\gamma'\in\n A(\Omega'_{\D})$, and consider the collection $\{\eta_{\gamma}\}$ such that $\eta_{\gamma}=0$ if $\gamma\neq\gamma'$ and $\eta_{\gamma'}=1$. Then (\ref{dsveq}) becomes
\begin{gather*}
\Big(b_n(v_{\gamma'}(k))\Big)_{\bar t}+\sum\limits_{i=1}^d\left(-\frac{v_{\gamma' x_i}(k)}h\right)+\sum\limits_{i\text{ s.t. }\gamma'-e_i\in\n A}\frac{v_{(\gamma'-e_i)x_i}(k)}h-f^{\D}_{(\gamma',k)}=0\\[4mm]\iff
\Big(b_n(v_{\gamma'}(k))\Big)_{\bar t}+\sum\limits_{i=1}^d\left(-\frac{v_{\gamma' x_i}(k)}h\right)+\sum\limits_{i=1}^d\frac{v_{(\gamma'-e_i)x_i}(k)}h-f^{\D}_{(\gamma',k)}=0
\end{gather*}
which is (\ref{dsvequiv}) for $\gamma'$. Since $\gamma'$ was arbitrary in $\n A(\Omega'_{\D})$, it follows (ii)' is satisfied.\hfill{$\square$}

\begin{lemma}\label{exuniqdsv} Fix a discretization $\D$ with small $h$. Then for any $R>0$, to each $[f]_{\D}\in\n F^R_{\D}$ there corresponds a unique discrete state vector.
\end{lemma}

\emph{Proof.} First we prove uniqueness. Let $[v([f_{\D}])]_{\D}, [\tilde v([f_{\D}])]_{\D}$ both satisfy Definition \ref{dsvdef}. At the outset it is clear that $v(0)=\tilde v(0)$ due to (i) and (iii). Proceeding by induction, fix $k,~1\leq k\leq n$ and suppose $v(k-1)=\tilde v(k-1)$. In (\ref{dsveq}) for both $v(k)$ and $\tilde v(k)$, plug in $\eta=v(k)-\tilde v(k)$, and subtract the resulting equalities. We obtain
\[
\sum\limits_{\n A}h^d\left[\Big((b_n(v_{\gamma}(k)))_{\bar t}-(b_n(\tilde v_{\gamma}(k)))_{\bar t}\Big)(v_{\gamma}(k)-\tilde v_{\gamma}(k))+\sum\limits_{i=1}^d\Big(v_{\gamma x_i}(k)-\tilde v_{\gamma x_i}(k)\Big)^2\right]=0,
\]
but we note that for each $\gamma$,
\begin{align*}
(b_n(v_{\gamma}(k)))_{\bar t}-(b_n(\tilde v_{\gamma}(k)))_{\bar t}&=\frac{b_n(v_{\gamma}(k))-b_n(v_{\gamma}(k-1))}{\tau}-\frac{b_n(\tilde v_{\gamma}(k))-b_n(\tilde v_{\gamma}(k-1))}{\tau}=\\[3mm]&=\frac{b_n(v_{\gamma}(k))-b_n(\tilde v_{\gamma}(k))}{\tau}
\end{align*}
per the induction hypothesis. It follows
\[
\sum\limits_{\n A}h^d\left[\frac{b_n(v_{\gamma}(k))-b_n(\tilde v_{\gamma}(k))}{\tau}(v_{\gamma}(k)-\tilde v_{\gamma}(k))+\sum\limits_{i=1}^d\Big(v_{\gamma x_i}(k)-\tilde v_{\gamma x_i}(k)\Big)^2\right]=0.
\]
Since $b$ is monotonically increasing, so is $b_n$. It follows that all terms in the above sum are non-negative, and so each term is identically $0$. In particular, due to the monotonicity of $b_n$ it follows that $v_{\gamma}(k)=\tilde v_{\gamma}(k)$ for $\gamma\in\n A(\Omega'_{\D})$. Due to (iii), this can be extended to $v(k)=\tilde v(k)$. By induction, this proves $v=\tilde v$.\\

Now we prove the existence. Fix a discretization $\D$ and $[f]_{\alpha}$. We will establish existence by induction on $k$. When $k=0$, we let $v(0)$ be given as in (i) and (iii) of Definition \ref{dsvdef}. By the induction hypothesis at level $k$, suppose that the first $k-$components $v(0),v(1),\ldots,v(k-1)$ have been constructed. We will give $v(k)$ now by the method of successive approximations. Obviously $v(k)$ on the lattice at the boundary of $\Omega_{\D}$ is just set to be $0$. For the lattice points in the interior, we  notice that (\ref{dsvequiv}) can be written in the following way
\begin{equation}\label{vsystem}
\frac1{\tau}\big[b_n(v_{\gamma}(k))-b_n(v_{\gamma}(k-1))\big]+\frac1{h^2}\left[2dv_{\gamma}(k)-\sum\limits_{i=1}^d(v_{\gamma+e_i}(k)+v_{\gamma-e_i}(k))\right]=f^{\D}_{(\gamma,k)}.
\end{equation}
So set $v^0=v(k-1)$, and having calculated $v^N$, obtain $v^{N+1}$ from the following system of equations:
\begin{equation}\label{vNsystem}
\frac1{\tau}b_n(v_{\gamma}^{N+1})+\frac{2d}{h^2}v_{\gamma}^{N+1}=\frac1{\tau}b_n(v_{\gamma}(k-1))+\frac1{h^2}\sum\limits_{i=1}^d(v_{\gamma+e_i}^N+v_{\gamma-e_i}^N)+f^{\D}_{(\gamma,k)}.
\end{equation}
Since the left hand side of \eqref{vNsystem} is monotonically increasing with respect to $v_{\gamma}^{N+1}$ and has a range $\mathbb{R}$, there is a unique solution 
$v^{N+1}$, and hence the sequence $\{v^N\}$ is well-defined. Now for each $\gamma$, subtract (\ref{vNsystem}) for $N$ and $N-1$ to get the system
\begin{equation}\label{system2}
\frac{h^2}{\tau}\Big(b_n(v_{\gamma}^{N+1})-b_n(v_{\gamma}^N)\Big)+2d\Big(v_{\gamma}^{N+1}-v_{\gamma}^N\Big)=\sum\limits_{i=1}^d\Big[(v_{\gamma+e_i}^N-v_{\gamma+e_i}^{N-1})+(v_{\gamma-e_i}^N-v_{\gamma-e_i}^{N-1})\Big].
\end{equation}
Now let
\begin{equation}\label{zeta}
\zeta_{\D,N}^{\gamma}:=\int\limits_0^1b_n'\Big(\theta v_{\gamma}^{N+1}+(1-\theta)v_{\gamma}^N\Big)\,d\theta,
\end{equation}
so it follows
\[
b_n(v_{\gamma}^{N+1})-b_n(v_{\gamma}^N)=\zeta_{\D,N}^{\gamma}\Big(v_{\gamma}^{N+1}-v_{\gamma}^N\Big),
\]
and
\begin{equation}\label{zetabig}
\zeta_{\D,N}^{\gamma}\geq\inf\limits_{x\in\bb R}b_n'(x)\geq\bar b,
\end{equation}
independently of $N,\gamma,\D$. Hence, system (\ref{system2}) can be written as
\begin{equation}\label{system3}
v_{\gamma}^{N+1}-v_{\gamma}^N=\frac1{2d+\frac{h^2}{\tau}\zeta_{\D,N}^{\gamma}}\sum\limits_{i=1}^d\Big[(v_{\gamma+e_i}^N-v_{\gamma+e_i}^{N-1})+(v_{\gamma-e_i}^N-v_{\gamma-e_i}^{N-1})\Big].
\end{equation}
By (\ref{zetabig}) we have that
\[
\frac1{2d+\frac{h^2}{\tau}\zeta_{\D,N}^{\gamma}}\leq\frac1{2d+\frac{h^2}{\tau}\bar b}
\]
uniformly over $\gamma$ and $N$. Let
\[
A_N:=\max\limits_{\gamma}|v_{\gamma}^{N+1}-v_{\gamma}^N|,
\]
then (\ref{system3}) implies that
\[
|v_{\gamma}^{N+1}-v_{\gamma}^N|\leq\frac{2d}{2d+\frac{h^2}{\tau}\bar b}A_{N-1}
\]
for each $\gamma$. Define
\[
\delta:=\frac{2d}{2d+\frac{h^2}{\tau}\bar b}.
\]
It is clear that $\delta\in(0,1)$. Thus we can arrive at the chain of inequalities
\begin{equation}\label{chain}
A_N\leq\delta A_{N-1}\leq\delta^2 A_{N-2}\leq\cdots\leq\delta^N A_0.
\end{equation}
Now, for any $N>M\geq0$, for fixed $\gamma$ we can write
\[
v_{\gamma}^N=v_{\gamma}^M+\sum\limits_{\ell=M}^{N-1}\big(v_{\gamma}^{\ell+1}-v_{\gamma}^{\ell}\big),
\]
which implies that
\begin{equation}\label{telescope}
|v_{\gamma}^N|\leq|v_{\gamma}^M|+\sum\limits_{\ell=M}^{N-1}|v_{\gamma}^{\ell+1}-v_{\gamma}^{\ell}|\leq|v_{\gamma}^M|+A_0\sum\limits_{\ell=M}^{N-1}\delta^{\ell}\leq|v_{\gamma}^M|+A_0\sum\limits_{\ell=M}^{\infty}\delta^{\ell}.
\end{equation}
Setting $M=0$ in (\ref{telescope}) gives that the sequence $\{v_{\gamma}^N\}$ is uniformly bounded in $\bb R$ with respect to $N$. Now let $\{v_{\gamma}^{N_p}\}$ be a subsequence which converges to $\liminf\limits_{N\ra\infty}v_{\gamma}^N$. Choose $M=N_p$ in an inequality similar to (\ref{telescope}) to see that
\begin{equation*}
v_{\gamma}^N\leq v_{\gamma}^{N_p}+A_0\sum\limits_{\ell=N_p} ^{\infty}\delta^{\ell},\qquad\forall N>N_p,
\end{equation*}
so that
\[
\limsup\limits_{N\ra\infty}v_{\gamma}^N\leq v_{\gamma}^{N_p}+A_0\sum\limits_{\ell=N_p} ^{\infty}\delta^{\ell},\quad p=1,2,\ldots
\]
which implies, upon sending $p\ra\infty$ that
\[
\limsup\limits_{N\ra\infty}v_{\gamma}^N\leq\liminf\limits_{N\ra\infty}v_{\gamma}^N
\]
and so the sequence $\{v_{\gamma}^N\}$ converges to a finite limit, for each $\gamma$. It follows we can define
\begin{equation}\label{successive}
v_{\gamma}(k)=\lim\limits_{N\ra\infty}v_{\gamma}^N,\quad\gamma\in\n A(\Omega_{\D}').
\end{equation}
We claim that $v(k)$ given by (\ref{successive}) satisfies (\ref{dsveq}). Due to Proposition \ref{equiv}, it is enough to see whether $v(k)$ satisfies system (\ref{vsystem}). But this follows immediately since $b_n$ and the identity map are continuous functions. This finishes the step of the induction, and therefore the proof.\hfill{$\square$}

The next lemma formulates the necessary and sufficient condition for the convergence of the discrete optimal control problems to the continuous optimal control problem.
\begin{lemma}\label{Vasil}\cite{Vasilev1} The sequence of discrete optimal control problems $\m I_{n}$ approximates the continuous optimal control problem $\m I$ with respect to the functional if and only if the following conditions are satisfied:
	\begin{description}
		\item{(i)} For any $f\in\n F^R$, we have $\n Q_{\D}(f)\in\n F_{\D}^R$, and 
		\begin{equation}\label{firstcondition}
			\limsup\limits_{\D\ra0}\Big( \n I_{\D}(\n{Q}_{\D}(f))-\n{J}(f) \Big ) \le 0.
		\end{equation}
		\item{(ii)} For any $[f]_{\D}\in\n F_{\D}^R$, we have $\n P_{\D}([f]_{\D})\in\n F^R$, and
		\begin{equation}\label{secondcondition}
			\limsup\limits_{\D\ra0} \Big ( \n{J}(\n{P}_{\D}([f]_{\D})) -\n{I}_{\D}([f]_{\D})  \Big ) \le 0.
		\end{equation}
		\end{description}
\end{lemma}
\begin{proposition}\label{PQ} The maps $\n P_{\D}$ and $\n Q_{\D}$ satisfy the conditions of Lemma \ref{Vasil}.
\end{proposition}

\emph{Proof.} Fix $\ep>0$ and $\D$ arbitrary. First let $f\in\n F^{R}$. Then we note
\[
\Vert\n Q_{\D}(f)\Vert_{\ell_{\infty}}=\max\limits_{\n A(\n C_{\D}^D)}|f_{\alpha}|=\max\limits_{\n A(\n C_{\D}^D)}\left|\frac1{\tau h^d}\int\limits_{C_{\D}^{\alpha}}f(x,t)\,dx\,dt\right|\leq\Vert f\Vert_{L_{\infty}(D)}\leq R.
\]
Now let $[f]_{\D}\in\n F_{\D}^R$. We see
\[
\Vert\n P_{\D}([f]_{\D})\Vert_{L_{\infty}(D)}=\esssup\limits_{ D}|f^{\D}(x,t)|=\max\limits_{\n A(\n C_{\D}^D)}|f_{\alpha}|=\Vert[f]_{\D}\Vert_{\ell_{\infty}}\leq R,
\]
which completes the proof.  \ \hfill{$\square$}

The following proposition is proved in \cite{LSU} for a wider class of solutions than that given in Definition \ref{weaksoldef}:
\begin{remark}\label{some} It is proved in \cite{LSU} that there exists a unique solution to the Stefan problem in the sense of Definition \ref{weaksoldef}. Moreover, it is proved that if a function $v\in\overset{\circ}{W}{}_2^{1,1}(D)\cap L_{\infty}(D)$ satisfies integral identity (\ref{weaksol}) for \emph{some} functions $B,B_0$ of type $\n B$ and any admissible test function $\psi$, then it follows that $v$ is the unique weak solution to the Stefan Problem in the sense of Definition \ref{weaksoldef}.
\end{remark}

\begin{proposition}\label{phistrong} For any $\ep>0$, there exists $\delta>0$ such that
\begin{equation}\label{phibound}
\sum\limits_{\n A}h^d\sum\limits_{i=1}^d|\Phi_{\gamma x_i}|^2\leq(1+\ep)\Vert D\Phi\Vert_{L_2(\Omega)}^2
\end{equation}
whenever $h<\delta$.
\end{proposition}

\emph{Proof.} Fix $i\in\{1,2,\ldots,d\}$. For each $h>0$, define the function $\tilde\Phi_h^i$ as
\[
\tilde\Phi_h^i\Big|_{R_{\D}^{\gamma}}=\Phi_{\gamma x_i},~~\forall\gamma\in\n A,\qquad\tilde\Phi_h^i\equiv0~~\text{elsewhere on }\Omega.
\]
We will prove that
\[
\tilde\Phi_h^i\longrightarrow\frac{\partial\Phi}{\partial x_i}~~\text{strongly in }L_2(\Omega)~~\text{as }h\ra0.
\]
As an element of $W_2^1(\Omega)$, almost all restrictions of $\Phi$ to lines parallel to the $x_i$ direction are absolutely continuous. Let $z=(z_1,\ldots,z_d)$, and write
\[
z+he_i=(z_1,\ldots,z_i+h,\ldots,z_d).
\]
Then we have that for almost every $z\in\Omega$,
\begin{equation}\label{abscont}
\Phi(z+he_i)-\Phi(z)=\int\limits_{z_i}^{z_i+h}\frac{\partial\Phi}{\partial x_i}(z_1,\ldots,y,\ldots,z_d)\,dy,
\end{equation}
and we will agree to write $(z\backslash z_i,y)$ in place of the vector $(z_1,\ldots,y,\ldots,z_d)$, to emphasize that the variable in the $i-th$ direction of the $z$ vector is replaced by $y$. Using the definition of the collection $\{\Phi_{\gamma}\}$, (\ref{abscont}), and the Cauchy-Schwartz inequality, we get
\begin{gather}
\left\Vert\tilde\Phi_h^i-\frac{\partial\Phi}{\partial x_i}\right\Vert_{L_2(\Omega)}^2=\int\limits_{\Omega_{\D}}\left|\tilde\Phi_{\D}^i-\frac{\partial\Phi}{\partial x_i}\right|^2dx~+\left\Vert\frac{\partial\Phi}{\partial x_i}\right\Vert_{L_2(\Omega\backslash\Omega_{\D})}^2=\nonumber\\[4mm]=\left\Vert\frac{\partial\Phi}{\partial x_i}\right\Vert_{L_2(\Omega\backslash\Omega_{\D})}^2+\sum\limits_{\n A}\int\limits_{R_{\D}^{\gamma}}\left|\Phi_{\gamma x_i}-\frac{\partial\Phi}{\partial x_i}(x)\right|^2dx=\nonumber\\[4mm]=\left\Vert\frac{\partial\Phi}{\partial x_i}\right\Vert_{L_2(\Omega\backslash\Omega_{\D})}^2+\sum\limits_{\n A}\int\limits_{R_{\D}^{\gamma}}\left|\frac1{h^{d+1}}\left[\int\limits_{R_{\D}^{\gamma+e_i}}\Phi(z)\,dz-\int\limits_{R_{\D}^{\gamma}}\Phi(z)\,dz\right]-\frac{\partial\Phi}{\partial x_i}(x)\right|^2dx=\nonumber\\[4mm]=\left\Vert\frac{\partial\Phi}{\partial x_i}\right\Vert_{L_2(\Omega\backslash\Omega_{\D})}^2+\sum\limits_{\n A}\int\limits_{R_{\D}^{\gamma}}\left|\frac1{h^{d+1}}\int\limits_{R_{\D}^{\gamma}}\big[\Phi(z+he_i)-\Phi(z)\big]\,dz-\frac{\partial\Phi}{\partial x_i}(x)\right|^2dx=\nonumber\\[4mm]=\left\Vert\frac{\partial\Phi}{\partial x_i}\right\Vert_{L_2(\Omega\backslash\Omega_{\D})}^2+\sum\limits_{\n A}\int\limits_{R_{\D}^{\gamma}}\frac1{h^{2d+2}}\left|\int\limits_{R_{\D}^{\gamma}}\left[\left(\int\limits_{z_i}^{z_i+h}\frac{\partial\Phi}{\partial x_i}(z\backslash z_i,y)\,dy\right)-h\frac{\partial\Phi}{\partial x_i}(x)\right]\,dz\right|^2dx\leq\nonumber\\[4mm]\leq\left\Vert\frac{\partial\Phi}{\partial x_i}\right\Vert_{L_2(\Omega\backslash\Omega_{\D})}^2+\sum\limits_{\n A}\int\limits_{R_{\D}^{\gamma}}\frac1{h^{d+1}}\int\limits_{R_{\D}^{\gamma}}\int\limits_{z_i}^{z_i+h}\left|\frac{\partial\Phi}{\partial x_i}(z\backslash z_i,y)-\frac{\partial\Phi}{\partial x_i}(x)\right|^2dy\,dz\,dx.\label{phicalc}
\end{gather}

Since $\partial\Phi/\partial x_i\in L_2(\Omega)$ and
\begin{equation}\label{measuredrops}
m_d(\Omega\backslash\Omega_{\D})\searrow0~~\text{as }h\searrow0,
\end{equation}
it follows by the absolute continuity of the integral that the first term on the right-hand side of (\ref{phicalc}) vanishes as $h\ra0$. Thus we focus on the second term. Recall that by $x_{\gamma}$ we denote the natural corner of the prism $R_{\D}^{\gamma}$. By an application of Fubini's Theorem we switch the order of the integration with respect to $y$ and $z_i$. Hence we observe
\begin{gather}
\sum\limits_{\n A}\int\limits_{R_{\D}^{\gamma}}\frac1{h^{d+1}}\int\limits_{R_{\D}^{\gamma}}\int\limits_{z_i}^{z_i+h}\left|\frac{\partial\Phi}{\partial x_i}(z\backslash z_i,y)-\frac{\partial\Phi}{\partial x_i}(x)\right|^2dy\,dz\,dx=\nonumber\\[4mm]=\sum\limits_{\n A}\int\limits_{R_{\D}^{\gamma}}\frac1{h^{d+1}}\left(\int\limits_{R_{\D}^{\gamma}}(z_i-x_{\gamma i})\,dz +\int\limits_{R_{\D}^{\gamma+e_i}}(x_{\gamma i}+2h-z_i)\,dz~~\left|\frac{\partial\Phi}{\partial x_i}(z)-\frac{\partial\Phi}{\partial x_i}(x)\right|^2\right)dx\leq\nonumber\\[4mm]\leq\frac1{h^d}\sum\limits_{\n A}\int\limits_{R_{\D}^{\gamma}}\left(\int\limits_{R_{\D}^{\gamma}}\,dz +\int\limits_{R_{\D}^{\gamma+e_i}}\,dz~~\left|\frac{\partial\Phi}{\partial x_i}(z)-\frac{\partial\Phi}{\partial x_i}(x)\right|^2\right)dx.\label{phicalc1}
\end{gather}

Now fix $\ep>0$. Since $C^1(\overline{\Omega+B_1(0)})$ is dense in $W_2^1(\Omega+B_1(0))$, it follows that we can choose a function $g\in C^1(\overline{\Omega+B_1(0)})$ depending on $\ep$ such that
\begin{equation}\label{gapprox}
\Vert\Phi-g\Vert_{W_2^1(\Omega+B_1(0))}^2<\frac1{12+6m_d(\Omega)}\ep.
\end{equation}
Add and subtract the terms $\frac{\partial g}{\partial x_i}(z),\frac{\partial g}{\partial x_i}(x)$ in the integrands to obtain that
\begin{gather*}
\frac1{h^d}\sum\limits_{\n A}\int\limits_{R_{\D}^{\gamma}}\left(\int\limits_{R_{\D}^{\gamma}}\,dz +\int\limits_{R_{\D}^{\gamma+e_i}}\,dz~~\left|\frac{\partial\Phi}{\partial x_i}(z)-\frac{\partial\Phi}{\partial x_i}(x)\right|^2\right)dx\leq\\[4mm]\leq I_1+I_2+I_3 \intertext{where} I_1=\frac3{h^d}\sum\limits_{\n A}\int\limits_{R_{\D}^{\gamma}}\left(\int\limits_{R_{\D}^{\gamma}}\,dz +\int\limits_{R_{\D}^{\gamma+e_i}}\,dz~~\left|\frac{\partial\Phi}{\partial x_i}(z)-\frac{\partial g}{\partial x_i}(z)\right|^2\right)dx,\\[4mm] I_2=\frac3{h^d}\sum\limits_{\n A}\int\limits_{R_{\D}^{\gamma}}\left(\int\limits_{R_{\D}^{\gamma}}\,dz +\int\limits_{R_{\D}^{\gamma+e_i}}\,dz~~\left|\frac{\partial g}{\partial x_i}(z)-\frac{\partial g}{\partial x_i}(x)\right|^2\right)dx,\\[4mm] I_3=\frac3{h^d}\sum\limits_{\n A}\int\limits_{R_{\D}^{\gamma}}\left(\int\limits_{R_{\D}^{\gamma}}\,dz +\int\limits_{R_{\D}^{\gamma+e_i}}\,dz~~\left|\frac{\partial g}{\partial x_i}(x)-\frac{\partial\Phi}{\partial x_i}(x)\right|^2\right)dx.
\end{gather*}

We estimate each of $I_1,I_2,I_3$. Since $g\in C^1(\overline{\Omega+B_1(0)})$, it follows that $\partial g/\partial x_i$ is uniformly continuous on $\Omega+B_1(0)$. Therefore, there exists $\delta=\delta(g,\ep)>0$ such that
\[
\left|\frac{\partial g}{\partial x_i}(z)-\frac{\partial g}{\partial x_i}(x)\right|^2<\frac1{12+6m_d(\Omega)}\ep
\]
whenever $|z-x|<\delta$. Let $h_{\ep}>0$ satisfy
\[
\sqrt{d+3}~~h_{\ep}<\delta.
\]
Then it follows that for each $h<h_{\ep}$, any $\gamma\in\n A$, and any $x,z\in R_{\D}^{\gamma}\cup R_{\D}^{\gamma+e_i}$,
\[
\left|\frac{\partial g}{\partial x_i}(z)-\frac{\partial g}{\partial x_i}(x)\right|^2<\frac1{12+6m_d(\Omega)}\ep.
\]
Therefore,
\[
I_1=\frac3{h^d}\sum\limits_{\n A}\left(\int\limits_{R_{\D}^{\gamma}}\,dz +\int\limits_{R_{\D}^{\gamma+e_i}}\,dz~~\left|\frac{\partial\Phi}{\partial x_i}(z)-\frac{\partial g}{\partial x_i}(z)\right|^2\right)\leq6\Vert\Phi-g\Vert_{W_2^1(\Omega+B_1(0))}^2,
\]
\[
I_2\leq\frac3{h^d}\sum\limits_{\n A}\int\limits_{R_{\D}^{\gamma}}\left(\int\limits_{R_{\D}^{\gamma}}\,dz +\int\limits_{R_{\D}^{\gamma+e_i}}\,dz~~\ep~\,dx\right)\leq\frac{m_d(\Omega)}{2+m_d(\Omega)}\ep,
\]
\[
I_3=6\sum\limits_{\n A}\int\limits_{R_{\D}^{\gamma}}\left|\frac{\partial\Phi}{\partial x_i}(x)-\frac{\partial g}{\partial x_i}(x)\right|^2\,dx\leq6\Vert\Phi-g\Vert_{W_2^1(\Omega)}^2.
\]
Due to (\ref{gapprox}), these calculations imply that 
\[
I_1+I_2+I_3<\ep,\qquad\forall h\leq h_{\ep}
\]
which shows that the left-hand side of (\ref{phicalc1}) drops to $0$ as $h\ra0$. This proves the strong convergence of $\tilde\Phi_h^i$ to $\partial\Phi/\partial x_i$ in $L_2(D)$. Since
\[
\Vert\tilde\Phi_h^i\Vert_{L_2(\Omega)}^2=\sum\limits_{\n A}\int\limits_{R_{\D}^{\gamma}}|\tilde\Phi_h^i(x)|^2\,dx=\sum\limits_{\n A}h^d|\Phi_{\gamma x_i}|^2,
\]
estimate (\ref{phibound}) follows after running the previous argument for each $i=1,\ldots,d$.\hfill{$\square$}

\section{Estimates}\label{estimates}

\begin{theorem}\label{boundedness} (Discrete Maximum Principle) For any $R>0$, any $[f]_{\D}\in\n F_{\D}^R$, and any $\D$, the discrete state vector $[v([f]_{\D})]_{\D}$ satisfies the following estimate:
\begin{equation}\label{boundedest}
\Vert[v]_{\D}\Vert_{\ell_{\infty}}:=\max\limits_{0\leq k\leq n}\max\limits_{\n A(\Omega_{\D})}|v_{\gamma}(k)|\leq e^T\max\left\{\frac1{\bar b}\Vert[f]_{\D}\Vert_{\ell_{\infty}}~,~\Vert\Phi\Vert_{L_{\infty}(\Omega)}\right\}
\end{equation}
\end{theorem}

\emph{Proof.} Fix a discretization $\D=(\tau,h)$ and $[f]_{\D}\in\n F_{\D}^R$. There corresponds the unique discrete state vector $[v([f]_{\D})]_{\D}$ by Lemma \ref{exuniqdsv}. Consider the following transformation of the discrete state vector:
\begin{equation}\label{udef}
u_{\gamma}(k):=v_{\gamma}(k)e^{-t_k},\quad\forall(\gamma,k)\in\n A(D_{\D}).
\end{equation}

Then (\ref{mvt}) gives
\begin{align*}
\frac{b_n(v_{\gamma}(k))-b_n(v_{\gamma}(k-1))}{\tau}&=\zeta_{\D,k}^{\gamma}\frac{u_{\gamma}(k)e^{t_k}-u_{\gamma}(k-1)e^{t_{k-1}}}{\tau}=\\[3mm]&=\zeta_{\D,k}^{\gamma}\Big(u_{\gamma}(k)\frac1{\tau}(e^{t_k}-e^{t_{k-1}})+u_{\gamma}(k)_{\bar t}\,e^{t_{k-1}}\Big)=\\[3mm]&=\zeta_{\D,k}^{\gamma}\Big(u_{\gamma}(k)e^{t^k}+u_{\gamma}(k)_{\bar t}\,e^{t_{k-1}}\Big)
\end{align*}
where $t^k\in[t_{k-1},t_k]$ satisfies
\[
e^{t_k}-e^{t_{k-1}}=e^{t^k}\tau,
\]
and such a $t^k$ exists for each $k$ due to the Mean Value Theorem. It follows (\ref{dsvequiv}) is transformed as
\begin{equation*}\label{dsvu}
\zeta_{\D,k}^{\gamma}u_{\gamma}(k)e^{t^k}+\zeta_{\D,k}^{\gamma}u_{\gamma}(k)_{\bar t}\,e^{t_{k-1}}-e^{t_k}\Delta_h u_{\gamma}(k)=f_{\gamma,k}^{\D}
\end{equation*}
which yields
\begin{equation}\label{dsvu2}
\zeta_{\D,k}^{\gamma}u_{\gamma}(k)e^{t^k}+\zeta_{\D,k}^{\gamma}u_{\gamma}(k)_{\bar t}\,e^{t_{k-1}}-e^{t_k}\sum\limits_{i=1}^d\frac1h\big(u_{\gamma}(k)_{x_i}-u_{\gamma-e_i}(k)_{x_i}\big)=f_{\gamma,k}^{\D}.
\end{equation}

Now, if $u_{\gamma}(k)\leq 0$ for every $\alpha\in\n A(D_{\D})$, then $\max\limits_{\n A(D_{\D})}u_{\gamma}(k)\leq0$. If, instead, we have $u_{\gamma}(k)>0$ for some $\alpha=(\gamma,k)\in\n A(D_{\D})$, then $\max\limits_{\n A(D_{\D})}u_{\gamma}(k)>0$, and let $\alpha^*=(\gamma^*,k^*)$ be such that
\[
u_{\gamma^*}(k^*)=\max\limits_{\n A(D_{\D})}u_{\gamma}(k).
\]

By assumption, $\alpha^*$ cannot lie on $\n A(S_{\D})$ (i.e. the lateral boundary of $D_{\D}$). If $\alpha^*$ lies on $\n A(\Omega_{\D}^0)$, then we clearly have
\[
u_{\gamma^*}(k^*)=\max\limits_{\n A(\Omega_{\D})}\Phi_{\gamma}\leq\Vert\Phi\Vert_{L_{\infty}(\Omega)}.
\]
The final possibility is that $\alpha^*$ lies on $\n A(D'_{\D})$. In this case, (\ref{dsvu2}) is satisfied at $\alpha^*$, and moreover we must have
\[
u_{\gamma^*}(k^*)_{\bar t}\geq0,\qquad u_{\gamma^*}(k^*)_{x_i}\leq0~\forall i,\qquad u_{\gamma^*-e_i}(k^*)_{x_i}\geq0~\forall i
\]
by definition of $\alpha^*$. Per our assumptions,
\begin{equation}\label{zetabig1}
\zeta_{\D,k}^{\gamma}\geq\bar b
\end{equation}
uniformly for $k,\gamma,\D$. Hence (\ref{dsvu2}) yields the inequality
\[
u_{\gamma^*}(k^*)\leq\frac1{\bar b}f_{\alpha^*}^{\D}e^{-t^k}\leq\frac1{\bar b}\Vert[f]_{\D}\Vert_{\ell_{\infty}}.
\]
The past observations imply that
\[
\max\limits_{\n A(D_{\D})}v_{\gamma}(k)\leq e^{T}\max\left\{\frac1{\bar b}\Vert[f]_{\D}\Vert_{\ell_{\infty}}~,~\Vert\Phi\Vert_{L_{\infty}(\Omega)}\right\}
\]
In a completely analogous fashion we are able to obtain a uniform lower bound:
\[
\min\limits_{\n A(D_{\D})}v_{\gamma}(k)\geq e^T\min\left\{-\frac1{\bar b}\Vert[f]_{\D}\Vert_{\ell_{\infty}}~,~-\Vert\Phi\Vert_{L_{\infty}(\Omega)}\right\},
\]
giving (\ref{boundedest}).\hfill{$\square$}

\begin{theorem}\label{energy}(Discrete $W_2^{1,1}$ Energy Estimate) For $[f]_{\D}\in\n F_{\D}^R$ and any $\D$, the discrete state vector $[v([f]_{\D})]_{\D}$ satisfies the following estimate:
	\begin{gather}
	\sum\limits_{k=1}^n\tau\sum\limits_{\n A}h^d (v_{\gamma\bar t}(k))^2+\max\limits_{1\leq k\leq n}\sum\limits_{\n A}h^d\sum\limits_{i=1}^d (v_{\gamma x_i}(k))^2+\nonumber\\[4mm]+\sum\limits_{k=1}^n\tau^2\sum\limits_{\n A}h^d\sum\limits_{i=1}^d \big(v_{\gamma x_i\bar t}(k)\big)^2\leq\frac2{\bar b\min\{1,\bar b\}}~\Vert f^{\D}\Vert_{L_2(D)}^2+4\Vert\Phi\Vert_{W_2^1(\Omega)}^2~~=:\n E([f]_{\D})\label{energyest}
	\end{gather}
\end{theorem}

\emph{Proof.} In (\ref{dsveq}), for each $k=1,2,\ldots,n$ choose $\eta:=2\tau v_{\gamma}(k)_{\bar t}$, and consider the identity
\[
2\tau v_{\gamma x_i}(k)\big(v_{\gamma}(k)_{\bar t}\big)_{x_i}=(v_{\gamma x_i}(k))^2-(v_{\gamma x_i}(k-1))^2+\tau^2\big(v_{\gamma}(k)_{x_i\bar t}\big)^2.
\]
Upon using (\ref{mvt}) on the first term of (\ref{dsveq}) with the aforementioned $\eta$, we readily observe
\begin{equation}\label{eq1}
\sum\limits_{\n A}\left[2\tau\zeta_{\D,k}^{\gamma}(v_{\gamma\bar t}(k))^2+\sum\limits_{i=1}^d \Big((v_{\gamma x_i}(k))^2-(v_{\gamma x_i}(k-1))^2+\tau^2\big(v_{\gamma}(k)_{x_i\bar t}\big)^2\Big)-2\tau f_{\gamma,k}^{\D}v_{\gamma}(k)_{\bar t}\right]=0
\end{equation}
for each $k=1,2,\ldots, n$. Due to (\ref{zetabig1}), we can use Cauchy's Inequality with $\ep=\bar b$ on the last term to obtain from (\ref{eq1}) the inequality
\begin{equation}\label{eq2}
\sum\limits_{\n A}\left[\tau \bar b(v_{\gamma\bar t}(k))^2+\sum\limits_{i=1}^d(v_{\gamma x_i}(k))^2-(v_{\gamma x_i}(k-1))^2+\tau^2\big(v_{\gamma}(k)_{x_i\bar t}\big)^2\right]\leq\frac1{\bar b}\sum\limits_{\n A}\tau(f_{\gamma,k}^{\D})^2
\end{equation}
true for each $k=1,2,\ldots, n$. Perform a summation of (\ref{eq1}) over $k=1,2,\ldots,q\leq n$. We see
\begin{gather}
\sum\limits_{k=1}^q\tau\sum\limits_{\n A}h^d (v_{\gamma\bar t}(k))^2+\sum\limits_{\n A}h^d\sum\limits_{i=1}^d (v_{\gamma x_i}(q))^2+\sum\limits_{k=1}^q\tau^2\sum\limits_{\n A}h^d\sum\limits_{i=1}^d \big(v_{\gamma}(k)_{x_i\bar t}\big)^2\leq\nonumber\\[4mm]\leq\frac1{\bar b\min\{1,\bar b\}}\sum\limits_{k=1}^q\tau\sum\limits_{\n A}h^d(f_{\gamma,k}^{\D})^2+\frac1{\min\{1,\bar b\}}\sum\limits_{\n A}h^d\sum\limits_{i=1}^d (\Phi_{\gamma x_i})^2\label{summed}
\end{gather}

We note by the Cauchy-Schwartz inequality that
\[
\sum\limits_{k=1}^q\tau\sum\limits_{\n A}h^d(f_{\gamma,k}^{\D})^2\leq\int\limits_{D_{\D}}(f^{\D})^2\,dx\,dt\leq\Vert{f^{\D}}\Vert_{L_2(D)}^2,
\]
and owing to Proposition \ref{phistrong},
\begin{equation}\label{phiestimate}
\sum\limits_{\n A}h^d\sum\limits_{i=1}^d (\Phi_{\gamma x_i})^2\leq {2\Vert\Phi\Vert_{W_2^1(\Omega)}^2}.
\end{equation}
With these observations, choosing $q=n$ in (\ref{summed}) and $q$ the maximizer for the second term on the left-hand side, we arrive at the desired estimate.\hfill{$\square$}

\section{Theorem on Interpolations of a Discrete State Vector}\label{Interpolations}

We describe a few useful ways in which we can interpolate the discrete state vectors to functions over $D$. Recall that a discrete state vector assigns a unique value $v_{\gamma}(k)$ to each point in the lattice $\n L(\Omega_{\D})$. In particular, we can identify each cell in $\Omega_{\D}$ by its natural corner, which is a point in the aforementioned lattice. The collection of natural corners is indexed by the set $\n A$.\\

By $\tilde V_{\D}:D\ra\bb R$, it is meant an interpolation of a discrete state vector $[v]_{\D}$ which assigns to the interior and top face of each cell in $D_{\D}$ the value at its natural corner. That is,
\begin{equation}\label{pwconstant}
\tilde V_{\D}\Big|_{C_{\D}^{(\gamma,k)}}=v_{\gamma}(k),\qquad\tilde V_{\D}\Big|_{R_{\D}^{\gamma,k}}=v_{\gamma}(k),\qquad\forall(\gamma,k)\in\n A(\n C_{\D}^D),
\end{equation}
and we let $\tilde V_{\D}$ be $0$ elsewhere in $D$ that it is not already defined. Now for each $i=1,2,\ldots,d$, define the function $\tilde V_{\D}^i:D\ra\bb R$ as
\begin{equation}\label{derivativestep}
\tilde V_{\D}^i\Big|_{C_{\D}^{(\gamma,k)}}=v_{\gamma x_i}(k),\qquad\forall(\gamma,k)\in\n A(\n C_{\D}^D)
\end{equation}
and $0$ elsewhere in $D$ where it is not already given by (\ref{derivativestep}). Intuitively, the $\tilde V_{\D}^i$ are step functions which assign to each cell in $D_{\D}$ the value of the forward spatial difference at the natural corner. Next, for fixed $k=0,1,\ldots,n$, we define $V_{\D}^k:\Omega\ra\bb R$ as a spatial interpolation of the discrete state vector which assigns to each point in the lattice $\n L(\Omega_{\D})$ the corresponding value $v_{\gamma}(k)$, is linear with respect to any spatial variable when all other spatial variables are fixed, and is extended as $0$ on $\Omega\backslash\Omega_{\D}$. This gives a unique interpolation, and we note $V_{\D}^k$ is a continuous function. Then we define the function $V_{\D}:D\ra\bb R$ as the piece-wise constant interpolation of the functions $V_{\D}^k$ onto time. That is,
\begin{equation}\label{pwlinearpwconstant}
V_{\D}(x,t)=V_{\D}^k(x),\qquad t\in(t_{k-1},t_k],~k=1,2,\ldots,n
\end{equation}
and $V_{\D}(x,0)=V_{\D}^0(x)$. Finally, we define the function $V_{\D}':D\ra\bb R$ as the piece-wise linear interpolation of $V_{\D}^k$ onto time. That is,
\begin{equation}\label{pwlinear}
V_{\D}'(x,t)=V_{\D}^{k-1}(x)+\big(V_{\D}^k(x)\big)_{\bar t}(t-t_{k-1}),\qquad t\in[t_{k-1},t_k],~k=1,2,\ldots,n.
\end{equation}

Now let us make a few remarks about the spatial functions $V_{\D}^k$. Fix a rectangular prism $R_{\D}^{\gamma}$. Such a prism has $2^d$ vertexes, which are the elements of $\n L(R_{\D}^{\gamma})$. By definition of $V_{\D}^k$, one can see that for each $x\in R_{\D}^{\gamma}$, the value $V_{\D}^k(x)$ is the weighted average (with respect to distance from the point $x$ to each vertex) of the values $v_{\gamma^*}(k)$ where $\gamma^*\in\n A(R_{\D}^{\gamma})$. Therefore $V^k_{\D}$ satisfies the following representation in each prism $R_{\D}^{\gamma}$:
\begin{equation}\label{weighted1}
V_{\D}^k(x)=\sum\limits_{\gamma^*\in\n A(R_{\D}^{\gamma})}w_{\gamma^*}(x)v_{\gamma^*}(k),\qquad x\in R_{\D}^{\gamma}
\end{equation}
where each weight function $w_{\gamma^*}:R_{\D}^{\gamma}\ra[0,1]$ is continuous, and moreover we have
\begin{equation}\label{weightssumup1}
\sum\limits_{\gamma^*\in\n A(R_{\D}^{\gamma})}w_{\gamma^*}(x)=1,\qquad x\in R_{\D}^{\gamma},
\end{equation}
and we remark that, even though parts of the boundary of each prism intersects other prisms, the representation (\ref{weighted1}) is satisfied regardless of the prism chosen. Given (\ref{weighted1}), it easily follows that
\begin{equation}\label{remark1}
|V_{\D}^k|\Big|_{R_{\D}^{\gamma}}\leq\max\limits_{\gamma^*\in\n A(R_{\D}^{\gamma})}|v_{\gamma^*}(k)|,
\end{equation}
from which it is readily deduced that
\begin{equation}\label{remarkl}
\Vert V_{\D}^k\Vert_{L_{\infty}(\Omega)}\leq\max\limits_{\n A(\Omega_{\D}')}|v_{\gamma}(k)|.
\end{equation}

Continuing with the same set-up, fix a direction $i\in\{1,2,\ldots,d\}$. There are $2^{d-1}$ one-dimensional faces (i.e. lines connecting the vertexes) in $R_{\D}^{\gamma}$ which run parallel to the $x_i$ direction. To each of these lines corresponds a space-difference
\[
v_{\gamma^*x_i}(k),\qquad\gamma^*\in\n A(R_{\D}^{\gamma},i):=\n A(R_{\D}^{\gamma})\cap\{\gamma^*_i=\gamma_i\}.
\]
Then by construction, for each $x\in R_{\D}^{\gamma}$, the value $\frac{\partial}{\partial x_i}V_{\D}^k(x)$ is the weighted average (with respect to the distance from the point $x$ to each appropriate line) of the values $v_{\gamma^*x_i}(k)$ where $\gamma^*\in\n A(R_{\D}^{\gamma},i)$. Therefore we have the following representation for $\frac{\partial}{\partial x_i}V_{\D}^k$ in each prism $R_{\D}^{\gamma}$:
\begin{equation}\label{weighted2}
\frac{\partial}{\partial x_i}V_{\D}^k(x)=\sum\limits_{\gamma^*\in\n A(R_{\D}^{\gamma},i)}W_{\gamma^*}(x)v_{\gamma^* x_i}(k),\qquad x\in R_{\D}^{\gamma}
\end{equation}
where the weight functions $W_{\gamma^*}:R_{\D}^{\gamma}\ra[0,1]$ are continuous and satisfy
\begin{equation}\label{weightssumup2}
\sum\limits_{\gamma^*\in\n A(R_{\D}^{\gamma},i)}W_{\gamma^*}(x)=1,\qquad x\in R_{\D}^{\gamma}.
\end{equation}
It follows that
\begin{equation}\label{remark2}
\left|\frac{\partial}{\partial x_i}V_{\D}^k\right|~\Bigg|_{R_{\D}^{\gamma}}\leq\max\limits_{\gamma^*\in\n A(R_{\D}^{\gamma},i)}|v_{\gamma^*x_i}(k)|,\qquad\forall\gamma\in\n A.
\end{equation}
Using (\ref{remark2}), we estimate
\begin{gather}
\int\limits_{\Omega}\left|\frac{\partial}{\partial x_i}V_{\D}^k\right|^2\,dx=\int\limits_{\Omega_{\D}}\left|\frac{\partial}{\partial x_i}V_{\D}^k\right|^2\,dx\leq\sum\limits_{\gamma\in\n A}h^d\max\limits_{\gamma^*\in\n A(R_{\D}^{\gamma},i)}|v_{\gamma^*x_i}(k)|^2.\label{remark3}
\end{gather}
Since each line connecting lattice points is shared by $2^{d-1}$ rectangular prisms, (\ref{remark3}) allows us to conclude
\begin{equation}\label{remark4}
\int\limits_{\Omega}\left|\frac{\partial}{\partial x_i}V_{\D}^k\right|^2\,dx\leq 2^{d-1}\sum\limits_{\n A}h^d|v_{\gamma x_i}(k)|^2.
\end{equation}

\begin{theorem}\label{interp} Let $\{[f]_{\D}\}$ be a sequence of discrete control vectors such that there exists $R>0$ for which $[f]_{\D}\in\n F_{\D}^R$ for each $\D$. The following statements hold:
\begin{enumerate}[(a)]
	\item The sequences $\{\tilde V_{\D}\},\{V_{\D}\}, \{V_{\D}'\}$ are uniformly bounded in $L_{\infty}(D)$.
	
	\item For each $i\in\{1,\ldots,d\}$, the sequences $\{\tilde V_{\D}^i\}, \{\partial V_{\D}/\partial x_i\},\{\partial V_{\D}'/\partial x_i\}$ are uniformly bounded in $L_2(D)$. Moreover, the sequence $\{\partial V_{\D}'/\partial t\}$ is uniformly bounded in $L_2(D)$.
	
	\item The sequence $\{V_{\D}-V_{\D}'\}$ converges strongly to $0$ in $L_2(D)$ as $\tau\ra0$.
	
	\item {For each $k=1,\ldots,n$, the sequence $\{V_{\D}^k-\tilde V_{\D}(\cdot,t_k)\}$ converges strongly to $0$ in $L_2(\Omega)$ as $h\ra0$.} Furthermore, the sequence $\{\tilde V_{\D}-V_{\D}\}$ converges strongly to $0$ in $L_2(D)$ as $h\ra0$.
	
	\item For each $i\in\{1,2,\ldots,d\}$, the sequence $\{\partial V_{\D}/\partial x_i-\partial V_{\D}'/\partial x_i\}$ converges strongly to $0$ in $L_2(D)$ as $\tau\ra0$.
	
	\item For each $i\in\{1,2,\ldots,d\}$, the sequence $\{\tilde V_{\D}^i-\partial V_{\D}/\partial x_i\}$ converges weakly to $0$ in $L_2(D)$ as $\D\ra0$.
\end{enumerate}
\end{theorem}

\emph{Proof.} Due to Theorem \ref{boundedness}, (\ref{remarkl}), and the fact that $\Vert[f]_{\D}\Vert_{\ell_{\infty}}\leq R$ for each $\D$, statement (a) follows immediately. Now we move to prove statement (b). Fix $i\in\{1,2,\ldots,d\}$. We have
\begin{equation}
\int\limits_0^T\int\limits_{\Omega}|\tilde V_{\D}^i|^2\,dx\,dt=\sum\limits_{k=1}^n\tau\sum\limits_{\n A}h^d|v_{\gamma x_i}(k)|^2\leq C\n E([f]_{\D}),\label{derivative1}
\end{equation}
whence it is known each sequence $\{\tilde V_{\D}^i\}$ is uniformly bounded in $L_2(D)$. Next, due to (\ref{remark4}) we note
\begin{align*}
\int\limits_0^T\int\limits_{\Omega}\left|\frac{\partial}{\partial x_i}V_{\D}\right|^2\,dx\,dt=\sum\limits_{k=1}^n\tau\int\limits_{\Omega}\left|\frac{\partial}{\partial x_i}V_{\D}^k\right|^2\,dx\leq 2^{d-1}T\max\limits_{1\leq k\leq n}\sum\limits_{\n A}h^d|v_{\gamma x_i}(k)|^2.
\end{align*}
Adding the above inequality over $i=1,2,\ldots,d$ and using (\ref{energyest}), we obtain
\begin{equation}
\Vert D_xV_{\D}\Vert_{L_2(D)}^2\leq2^{d-1}T\max\limits_{1\leq k\leq n}\sum\limits_{\n A}h^d\sum\limits_{i=1}^d|v_{\gamma x_i}(k)|^2\leq C\n E([f]_{\D}),\label{derivative2}
\end{equation}
where $C$ is independent of $\D$. Now fix $i\in\{1,2,\ldots,d\}$ again. We observe
\begin{gather}
\int\limits_0^T\int\limits_{\Omega}\left|\frac{\partial}{\partial x_i}V_{\D}'\right|^2\,dx\,dt=\sum\limits_{k=1}^n~\int\limits_{t_{k-1}}^{t_k}\int\limits_{\Omega}\left|\frac{t_k-t}{\tau}\frac{\partial}{\partial x_i}V_{\D}^{k-1}(x)+\frac{t-t_{k-1}}{\tau}\frac{\partial}{\partial x_i}V_{\D}^k(x)\right|^2\,dx\leq\nonumber\\[4mm]\leq2\sum\limits_{k=1}^n\tau\int\limits_{\Omega}\left[\left|\frac{\partial}{\partial x_i}V_{\D}^{k-1}(x)\right|^2+\left|\frac{\partial}{\partial x_i}V_{\D}^{k}(x)\right|^2\right]\,dx\leq\nonumber\\[4mm]\leq2^{d+1}T\max\limits_{0\leq k\leq n}\sum\limits_{\n A}h^d|v_{\gamma}x_i(k)|^2.
\end{gather}
Adding the above inequality over $i=1,2,\ldots,d$ and recalling that $v_{\gamma}(0)=\Phi_{\gamma}$ for each $\gamma$, we arrive at
\begin{equation}
\Vert D_xV_{\D}'\Vert_{L_2(D)}^2\leq2^{d+1}T\max\limits_{0\leq k\leq n}\sum\limits_{\n A}h^d\sum\limits_{i=1}^d|v_{\gamma}x_i(k)|^2 \leq C\n E([f]_{\D}).\label{pwlinearxok}
\end{equation}
Now note that for each $k=1,\ldots, n$ and each $\gamma\in\n A$, we have due to (\ref{weighted1}) and the Cauchy-Schwartz inequality that
\begin{gather*}
|V_{\D\bar t}^k(x)|^2=\left|\frac{V_{\D}^k(x)-V_{\D}^{k-1}(x)}{\tau} \right|^2=\left|\frac1{\tau}\sum\limits_{\gamma^*\in\n A(R_{\D}^{\gamma})}w_{\gamma^*}(x)\Big(v_{\gamma^*}(k)-v_{\gamma^*}(k-1)\Big)\right|^2\leq\\[4mm]\leq\frac1{\tau^2}\left(\sum\limits_{\gamma^*\in\n A(R_{\D}^{\gamma})}w_{\gamma^*}^2(x)\right)\sum\limits_{\gamma^*\in\n A(R_{\D}^{\gamma})}|v_{\gamma^*}(k)-v_{\gamma^*}(k-1)|^2\leq\\[4mm]\leq\sum\limits_{\gamma^*\in\n A(R_{\D}^{\gamma})}|v_{\gamma^*\bar t}(k)|^2,\qquad\text{a.e. } x\in R_{\D}^{\gamma},
\end{gather*}
which allows us to deduce
\begin{gather}
\int\limits_0^T\int\limits_{\Omega}\left|\frac{\partial}{\partial t}V_{\D}'\right|^2\,dx\,dt=\sum\limits_{k=1}^n~\int\limits_{t_{k-1}}^{t_k}\int\limits_{\Omega_{\D}}|V_{\D\bar t}^k(x)|^2\,dx=\sum\limits_{k=1}^n\tau\sum\limits_{\n A}\int\limits_{R_{\D}^{\gamma}}|V_{\D\bar t}^k(x)|^2\,dx\leq\nonumber\\[4mm]\leq\sum\limits_{k=1}^n\tau\sum\limits_{\n A}\int\limits_{R_{\D}^{\gamma}}\sum\limits_{\gamma^*\in\n A(R_{\D}^{\gamma})}|v_{\gamma^*\bar t}(k)|^2\,dx\leq2^d\sum\limits_{k=1}^n\tau\sum\limits_{\n A}h^d|v_{\gamma\bar t}(k)|^2\,dx\label{eqn1}
\end{gather}
where the last inequality holds since each value $|v_{\gamma\bar t}(k)|^2$ for $\gamma\in\n A$ is summed up at most $2^d$ times in the $\n A(R_{\D}^{\gamma})$ summation (because each interior lattice point is shared by $2^d$ prisms). Thanks to the energy estimate (\ref{energyest}), it is then clear from (\ref{eqn1}) that
\begin{equation}\label{pwlineartok}
\left\Vert\frac{\partial}{\partial t}V_{\D}'\right\Vert_{L_2(D)}^2\leq2^d\n E([f]_{\D}),
\end{equation}
so ends the proof of statement (b).\\

Next we prove (c). To this end, note that for each $k=1,2,\ldots, n$ and $\gamma\in\n A$, we have
\begin{gather*}
|V_{\D}(x,t)-V_{\D}'(x,t)|^2=|V_{\D}^k(x)-V_{\D}^{k-1}(x)-V_{\D\bar t}^k(x)(t-t_{k-1})|^2=\\[4mm]=\left|\frac{t_k-t}{\tau}(V_{\D}^k(x)-V_{\D}^{k-1}(x))\right|^2\leq|V_{\D}^k(x)-V_{\D}^{k-1}(x)|^2=\\[4mm]=\left|\sum\limits_{\gamma^*\in\n A(R_{\D}^{\gamma})}w_{\gamma^*}(x)(v_{\gamma^*}(k) -v_{\gamma^*}(k-1))\right|^2\leq\sum\limits_{\gamma^*\in\n A(R_{\D}^{\gamma})}\tau^2|v_{\gamma^*\bar t}(k)|^2,\qquad\text{a.e. }(x,t)\in C_{\D}^{(\gamma,k)},
\end{gather*}
so that we can deduce
\begin{gather*}
\Vert V_{\D}-V_{\D}'\Vert_{L_2(D)}^2=\sum\limits_{(\gamma,k)\in\n A(\n C_{\D}^D)~}\int\limits_{C_{\D}^{(\gamma,k)}}|V_{\D}(x,t)-V_{\D}'(x,t)|^2\,dy\leq \\[4mm]\leq\sum\limits_{(\gamma,k)\in\n A(\n C_{\D}^D)~}\int\limits_{C_{\D}^{(\gamma,k)}}\sum\limits_{\gamma^*\in\n A(R_{\D}^{\gamma})}\tau^2|v_{\gamma^*\bar t}(k)|^2\,dy\leq\\[4mm]\leq2^d\tau^2\sum\limits_{k=1}^n\tau\sum\limits_{\n A}h^d|v_{\gamma\bar t}(k)|^2\longrightarrow0~~\text{as}~~\tau\ra0
\end{gather*}
thanks to Theorem \ref{energy}. This ends the proof of (c).\\

The proof of statement (d) follows. For each $k=1,2,\ldots,n$ and $\gamma\in\n A$, we observe
\begin{gather}
|\tilde V_{\D}(x,t_k)-V_{\D}^k(x)|^2=|v_{\gamma}(k)-V_{\D}^k(x)|^2=\left|v_{\gamma}(k)-\sum\limits_{\gamma^*\in\n A(R_{\D}^{\gamma})}w_{\gamma^*}(x)v_{\gamma^*}(k)\right|^2=\nonumber\\[4mm]=\left|\sum\limits_{\gamma^*\in\n A(R_{\D}^{\gamma})}w_{\gamma^*}(x)\big(v_{\gamma}(k)-v_{\gamma^*}(k)\big)\right|^2\leq\sum\limits_{\gamma^*\in\n A(R_{\D}^{\gamma})}\big|v_{\gamma}(k)-v_{\gamma^*}(k)\big|^2,\qquad\text{a.e. }x\in R_{\D}^{\gamma}\label{sum1}.
\end{gather}
We note that if $\gamma=(k_1,k_2,\ldots,k_d)$, then each $\gamma^*\in\n A(R_{\D}^{\gamma})$ satisfies that $\gamma^*_i\in\{k_i,k_i+1\}$. Therefore, for each fixed $\gamma^*\in\n A(R_{\D}^{\gamma})$, there is a (not necessarily unique) path along the edges of the prism $R_{\D}^{\gamma}$ which starts at $x_{\gamma}$, ends at $x_{\gamma^*}$, and is made up of gluing together at most $d$ one-dimensional edges of the prism. Call such a path $P_{\gamma\ra\gamma^*}$, and $T_P(x)$ the tangent vector to the path at point $x$. It is easy to see then that we can write
\begin{equation}\label{path}
v_{\gamma}(k)-v_{\gamma^*}(k)=\int\limits_{P_{\gamma\ra\gamma^*}}D_xV_{\D}\cdot T_P(x)\,dP=\sum\limits_{P_{\gamma\ra\gamma^*}} hv_{\gamma'x_j}(k)
\end{equation}
where the sum on the right-hand side of (\ref{path}) is taken over the $\gamma'$ that correspond to vertexes of $R_{\D}^{\gamma}$ which lie on the path $P_{\gamma\ra\gamma^*}$ (except for the end-point $x_{\gamma^*}$), and $j$ corresponds to the spatial direction that the path $P_{\gamma\ra\gamma^*}$ takes in moving from $x_{\gamma'}$ to the next vertex that lies on the path. With this observation in hand and using the Cauchy-Schwartz inequality, the following estimate is true, uniformly over the path chosen, and uniformly over $\gamma^*$:
\begin{equation}\label{nopath}
|v_{\gamma}(k)-v_{\gamma^*}(k)|^2\leq d\sum\limits_{\text{edges of }R_{\D}^{\gamma}} h^2|v_{\gamma'x_j(k)}|^2
\end{equation}
where the sum on the right-hand side of (\ref{nopath}) is taken over all $\gamma'$ and $j$ such that $\gamma'\in\n A(R_{\D}^{\gamma})$ and $\gamma'+e_j\in\n A(R_{\D}^{\gamma})$ (intuitively, recall that the spatial differences $v_{\gamma'x_j}$ are in natural bijection with the edges of the lattice. So effectively, the sum is over all edges of the prism $R_{\D}^{\gamma}$). Therefore, using (\ref{nopath}) and (\ref{sum1}), we have for each $k=1,\ldots,n$,
\begin{align}
|\tilde V_{\D}(x,t_k)-V_{\D}^k(x)|^2&\leq\sum\limits_{\gamma^*\in\n A(R_{\D}^{\gamma})} d\sum\limits_{\text{edges of }R_{\D}^{\gamma}} h^2|v_{\gamma'x_j}(k)|^2\nonumber\\[4mm]&\leq(2^d-1)d\sum\limits_{\text{edges of }R_{\D}^{\gamma}} h^2|v_{\gamma'x_j}(k)|^2,\qquad\text{a.e. }x\in R_{\D}^{\gamma}\label{edgesest}
\end{align}
since there are $2^d-1$ vertexes $x_{\gamma^*}$ other than $x_{\gamma}$ in $R_{\D}^{\gamma}$. By using (\ref{edgesest}) we derive
\begin{align*}
\Vert\tilde V_{\D}(\cdot,t_k)-V_{\D}^k\Vert_{L_2(\Omega)}^2&=\sum\limits_{\n A}\int\limits_{R_{\D}^{\gamma}}|\tilde V_{\D}(x,t_k)-V_{\D}'(x)|^2\,dx\,dt\leq\\[4mm]&\leq\sum\limits_{\n A}h^d(2^d-1)d\sum\limits_{\text{edges of }R_{\D}^{\gamma}} h^2|v_{\gamma'x_j}(k)|^2\leq\\[4mm]&\leq\sum\limits_{\n A}h^d(2^d-1)~d~2^{d-1}\sum\limits_{i=1}^d h^2|v_{\gamma x_i}(k)|^2,
\end{align*}
where the last inequality holds since each edge in the lattice is shared by at most $2^{d-1}$ prisms. Finally we deduce
\[
\Vert\tilde V_{\D}(\cdot,t_k)-V_{\D}^k\Vert_{L_2(\Omega)}^2\leq h^2 d(2^d-1)2^{d-1}\max\limits_{1\leq k\leq n}\sum\limits_{\n A}h^d\sum\limits_{i=1}^d|v_{\gamma x_i}(k)|^2\longrightarrow0~~\text{as}~~h\ra0,
\]
uniformly over $k$, where again we have made use of Theorem \ref{energy}. Since
\[
\Vert\tilde V_{\D}-V_{\D}\Vert_{L_2(D)}^2\leq T\max\limits_{1\leq k\leq n}\Vert\tilde V_{\D}(\cdot,t_k)-V_{\D}^k\Vert_{L_2(\Omega)}^2 \longrightarrow0~~\text{as}~~h\ra0,
\]
statement (d) follows.\\

Now we move to proving (e). In this regard, it will be enough to estimate $\Vert D_xV_{\D}-D_xV_{\D}'\Vert_{L_2(D)}^2$. So first fix $i=1,2,\ldots,d$. For each $k=1,2,\ldots,n$ and $\gamma\in\n A$, we see that for almost every $(x,t)\in C_{\D}^{(\gamma,k)}$,
\begin{gather}
\left|\frac{\partial}{\partial x_i}V_{\D}(x,t)-\frac{\partial}{\partial x_i}V_{\D}'(x,t)\right|^2=\nonumber\\[4mm]=\left|\frac{\partial}{\partial x_i}V_{\D}^k(x)-\frac{t_k-t}{\tau}\frac{\partial}{\partial x_i}V_{\D}^{k-1}(x)-\frac{t-t_{k-1}}{\tau}\frac{\partial}{\partial x_i}V_{\D}^k(x)\right|^2\leq\nonumber\\[4mm]\leq\left|\frac{\partial}{\partial x_i}V_{\D}^k(x)-\frac{\partial}{\partial x_i}V_{\D}^{k-1}(x)\right|^2=\nonumber\\[4mm]=\left|\sum\limits_{\gamma^*\in\n A(R_{\D}^{\gamma},i)}W_{\gamma^*}(x)v_{\gamma^* x_i}(k)-\sum\limits_{\gamma^*\in\n A(R_{\D}^{\gamma},i)}W_{\gamma^*}(x)v_{\gamma^* x_i}(k-1)\right|^2\leq\nonumber\\[4mm]\leq\sum\limits_{\gamma^*\in\n A(R_{\D}^{\gamma},i)}|v_{\gamma^* x_i}(k)-v_{\gamma^* x_i}(k-1)|^2=\sum\limits_{\gamma^*\in\n A(R_{\D}^{\gamma},i)}\tau^2|v_{\gamma^* x_i\bar t}(k)|^2.\label{sum2}
\end{gather}
Hence,
\begin{gather*}
\left\Vert\frac{\partial}{\partial x_i} V_{\D}-\frac{\partial}{\partial x_i}V_{\D}'\right\Vert_{L_2(D)}^2=\sum\limits_{k=1}^n~\int\limits_{t_{k-1}}^{t_k}\sum\limits_{\n A}\int\limits_{R_{\D}^{\gamma}}\left|\frac{\partial}{\partial x_i}V_{\D}(x,t)-\frac{\partial}{\partial x_i}V_{\D}'(x,t)\right|^2\,dx\,dt\leq\\[4mm]\leq\sum\limits_{k=1}^n\tau\sum\limits_{\n A}h^d\sum\limits_{\gamma^*\in\n A(R_{\D}^{\gamma},i)}\tau^2|v_{\gamma^* x_i\bar t}(k)|^2\leq\\[4mm]\leq\sum\limits_{k=1}^n\tau2^{d-1}\sum\limits_{\n A}h^d\tau^2|v_{\gamma x_i\bar t}(k)|^2,
\end{gather*}
where the last inequality holds since each edge of the prism $R_{\D}^{\gamma}$ is shared by at most $2^{d-1}$ prisms. Thus,
\[
\left\Vert D_xV_{\D}-D_xV_{\D}'\right\Vert_{L_2(D)}^2\leq \tau2^{d-1}\sum\limits_{k=1}^n\tau^2\sum\limits_{\n A}h^d\sum\limits_{i=1}^d|v_{\gamma x_i\bar t}(k)|^2\longrightarrow0~~\text{as}~~\tau\ra0
\]
due to Theorem \ref{energy}. Statement (e) follows.\\

Moving on to statement (f), fix $i\in\{1,2,\ldots,d\}$. We will now prove that the sequence $\{\tilde V_{\D}^i-\partial V_{\D}/\partial x_i\}$ converges weakly to $0$ in $L_2(D)$. Due to (b), it is clear that both sequences $\{\tilde V_{\D}^i\},\{\partial V_{\D}/\partial x_i\}$ have weak limit points in $L_2(D)$. So let $g,g^*\in L_2(D)$ be weak limit points of $\{\tilde V_{\D}^i\}$,$\{\partial V_{\D}/\partial x_i\}$ in $L_2(D)$ respectively. In particular, $\tilde V_{\D'}^i\ra g$ weakly in $L_2(D)$ as $\D'\ra0$, where $\D'$ is some subsequence of $\D$.  
Let us fix a step-function on $D$, which is of the form
\begin{equation}\label{step}
s(y)=\sum\limits_{j=1}^ma_j\chi_{E_j}(y)
\end{equation}
where $E_j$'s are formed with intersections of $D$ with rectangles in $\bb R^{d+1}$, $E_j$'s partition $D$, $\chi_E$ is the characteristic function of the set $E\subset D$, and $a_j\in\bb R$ for each $j=1,\ldots,m$. Since the class of such step functions is dense in $L_2(D)$ it is satisfactory to prove the claim (f) for arbitrary step function $s$ of type \eqref{step}. Recall that here and in the sequel, $y=(x,t)$ and $\alpha=(\gamma,k)$. Since $m_{d+1}(\partial D)=0$ (where $m_{d+1}$ is the Lebesgue measure on $\bb R^{d+1}$), it follows by this construction that $m_{d+1}(\partial E_j)=0$ for each $j=1,\ldots,m$, and therefore the set
\[
\partial E:=\bigcup_{j=1}^m\partial E_j
\]
has $d+1$-st dimensional Lebesgue measure $0$. For each $k=1,2,\ldots,n$ and $\gamma\in\n A$, we observe
\begin{align}
\tilde V_{\D}^i(x,t)-\frac{\partial}{\partial x_i}V_{\D}(x,t)&=v_{\gamma x_i}(k)-\sum\limits_{\gamma^*\in\n A(R_{\D}^{\gamma},i)}W_{\gamma^*}(x)v_{\gamma^* x_i}(k)=\nonumber\\[4mm]&=\sum\limits_{\gamma^*\in\n A(R_{\D}^{\gamma},i)}W_{\gamma^*}(x)\big(v_{\gamma x_i}(k)-v_{\gamma^* x_i}(k)\big)=\nonumber\\[4mm]&=\sum\limits_{\gamma^*\in\n A(R_{\D}^{\gamma},i)}W_{\gamma^*}(x)\Big(\tilde V_{\D}^i(x,t)-\tilde V_{\D}^i(x+hz_{\gamma^*},t)\Big),\qquad\forall(x,t)\in C_{\D}^{(\gamma,k)},\nonumber
\end{align}
where $z_{\gamma^*}:=\frac1h\big(x_{\gamma^*}-x_{\gamma}\big)$. Therefore,
\begin{gather}
\int\limits_0^T\int\limits_{\Omega}\left(\tilde V_{\D}^i(x,t)-\frac{\partial}{\partial x_i}V_{\D}(x,t)\right)s(x,t)\,dx\,dt=\nonumber\\[4mm]=\sum\limits_{k=1}^n\int\limits_{t_{k-1}}^{t_k}\sum\limits_{\n A}\sum\limits_{\gamma^*\in\n A(R_{\D}^{\gamma},i)}\int\limits_{R_{\D}^{\gamma}}W_{\gamma^*}(x)\Big(v_{\gamma x_i}(k)-v_{\gamma^* x_i}(k)\Big)\sum\limits_{j=1}^ma_j \chi_{E_j}(x,t)\,dx\,dt\label{great}.
\end{gather}
We now intend to switch the order of the summations on the right-hand side of (\ref{great}). To do this, recall that the summation over $\gamma^*$ is taken over all indexes that correspond to vertexes of the prism $R_{\D}^{\gamma}$ which satisfy $\gamma^*_i=\gamma_i$. Since all prisms $R_{\D}^{\gamma}$ are congruent, it follows the vector $hz_{\gamma^*}$ that connects $x_{\gamma}$ to $x_{\gamma^*}$ does not depend on the specific coordinates of $\gamma$ or $\gamma^*$; it only depends on their difference (which is itself independent of $\D$). Since $|\n A(R_{\D}^{\gamma,i})|=2^{d-1}$, the vectors $z_{\gamma^*}$ are taken from the set
\[
\n Z:=\Big\{z\in\{0,1\}^d~\Big|~i\text{-th coordinate of }z\text{ is equal to }0\Big\}.
\]
Consequently, the summation over $\gamma^*\in\n A(R_{\D}^{\gamma},i)$ can be thought of as a summation over the elements of $\n Z$, since $\n Z$ is in bijection with $\n A(R_{\D}^{\gamma},i)$. Let $\gamma^z$ be the unique index in $\n A(R_{\D}^{\gamma},i)$ that is identified by $z$. We remark that the set $\n Z$ is independent of $\D$. Moreover, we can identify $W_{\gamma^*}(x)$ purely by the corresponding $z\in\n Z$, so we write $W_z(x):=W_{\gamma^z}(x)$.\\

It follows from (\ref{great}) that
\begin{gather}
\int\limits_0^T\int\limits_{\Omega}\left(\tilde V_{\D}^i(x,t)-\frac{\partial}{\partial x_i}V_{\D}(x,t)\right)s(x,t)\,dx\,dt=\nonumber\\[4mm]=\sum\limits_{z\in\n Z}\sum\limits_{(\gamma,k)\in\n A(\n C_{\D}^D)}\Big(v_{\gamma x_i}(k)-v_{\gamma^z x_i}(k)\Big)\sum\limits_{j=1}^ma_j\int\limits_{C_{\D}^{(\gamma,k)}\cap E_j}W_z(x)\,dy\label{greater}.
\end{gather}

Now fix $z\in\n Z$. Define
\[
C_z:=\frac1{\tau h^d}\int\limits_{C_{\D}^{\alpha}}W_z(x)\,dy
\]
and, since $W_z(x)$ is non-negative and either linear or constant in each variable $x_1,\ldots,x_d$, it follows that $C_z\in[0,1]$, and in particular $C_z$ is independent of $\alpha$. Define the set
\[
\n C_{s}:=\Big\{C_{\D}^{\alpha}\in\n C_{\D}^D~\Big|~~\forall j=1,\ldots,m,~~m_{d+1}(C_{\D}^{\alpha}\backslash E_j)\neq0\Big\}.
\]
Intuitively, $\n C_s$ is the set of all cells in $\n C_{\D}^D$ whose interiors are not contained in a single $E_j$. Define $D_{\n C_s}:=\bigcup_{\n C_s}C_{\D}^{\alpha}$. Moreover, to each cell in $\n C_{\D}^D\backslash\n C_s$ we specify by $j_{\alpha}$ the unique index for which the interior of $C_{\D}^{\alpha}$ is contained in $E_{j_{\alpha}}$. Thus it is seen that
\[
\sum\limits_{j=1}^ma_j\int\limits_{C_{\D}^{\alpha}\cap E_j}W_z(x)\,dy=a_{j_{\alpha}}C_z\tau h^d=C_z\int\limits_{C_{\D}^{\alpha}}a_{j_{\alpha}}\,dy=C_z\int\limits_{C_{\D}^{\alpha}}s(y)\,dy,\qquad\forall\alpha\in\n A(\n C_{\D}^D\backslash\n C_s).
\]
We can write
\begin{align}
&\sum\limits_{(\gamma,k)\in\n A(\n C_{\D}^D)}\Big(v_{\gamma x_i}(k)-v_{\gamma^z x_i}(k)\Big)\sum\limits_{j=1}^ma_j\int\limits_{C_{\D}^{(\gamma,k)}\cap E_j}W_z(x)\,dy=\nonumber\\[4mm]&=\sum\limits_{\n A(\n C_{\D}^D\backslash\n C_s)}\Big(v_{\gamma x_i}(k)-v_{\gamma^z x_i}(k)\Big)\sum\limits_{j=1}^ma_j\int\limits_{C_{\D}^{\alpha}\cap E_j}W_z(x)\,dy+\nonumber\\[4mm]&+\sum\limits_{\n A(\n C_s)}\Big(v_{\gamma x_i}(k)-v_{\gamma^z x_i}(k)\Big)\sum\limits_{j=1}^ma_j\int\limits_{C_{\D}^{\alpha}\cap E_j}\big(W_z(x)-C_z+C_z)\,dy=\nonumber\\[4mm]&=\sum\limits_{\n A(\n C_{\D}^D)}\Big(v_{\gamma x_i}(k)-v_{\gamma^z x_i}(k)\Big)C_z\int\limits_{C_{\D}^{(\gamma,k)}}s(y)\,dy+I_{\D}=\nonumber\\[4mm]&=C_z\int\limits_0^T\int\limits_{\Omega}\Big(\tilde V_{\D}^i(x,t)-\tilde V_{\D}^i(x+hz,t)\Big)s(x,t)\,dx\,dt+I_{\D}\label{calc}
\end{align}
where
\[
I_{\D}:=\sum\limits_{\n A(\n C_s)}\Big(v_{\gamma x_i}(k)-v_{\gamma^z x_i}(k)\Big)\int\limits_{C_{\D}^{\alpha}}\big(W_z(x)-C_z)s(y)\,dy.
\]

It can be shown that $|I_{\D}|\ra0$ as $\D\ra0$. To see this, use the Cauchy-Schwartz inequality and Theorem \ref{energyest} to get 
\begin{align}
|I_{\D}|^2&\leq\left(\sum\limits_{\n A(\n C_s)}\Big|v_{\gamma x_i}(k)-v_{\gamma^z x_i}(k)\Big|^2\right)\sum\limits_{\n A(\n C_s)}\left(\int\limits_{C_{\D}^{\alpha}}\big(W_z(x)-C_z)s(y)\,dy\right)^2\leq\nonumber\\[4mm]&\leq\left(4\sum\limits_{\n A(\n C_{\D}^D)}|v_{\gamma x_i}(k)|^2\right)~\sum\limits_{\n A(\n C_s)}\tau h^d\int\limits_{C_{\D}^{\alpha}}\big|W_z(x)-C_z|^2|s(y)|^2\,dy\leq\nonumber\\[4mm]&\leq4\left(\sum\limits_{k=1}^n\tau\sum\limits_{\n A}h^d|v_{\gamma x_i}(k)|^2\right)4\sum\limits_{\n A(\n C_s)}\int\limits_{C_{\D}^{\alpha}}|s(y)|^2\,dy\leq\nonumber\\[4mm]&\leq16\n E([f]_{\D})\Vert s\Vert_{L_2(D_{\n C_s})}^2\leq C\Vert s\Vert_{L_2(D_{\n C_s})}^2.\label{calc1}
\end{align}
We claim that $m_{d+1}(D_{\n C_s})\ra0$ as $\D\ra0$. Consider the sets
\[
\partial E_{\delta}:=\Big\{y\in D~\big|~\text{dist}(y,\partial E)<\delta\},
\]
which are open in $\bb R^{d+1}$. Then $m_{d+1}(\partial E_{\delta})\searrow0$ as $\delta\searrow0$, since $m_{d+1}(\partial E)=0$. Now fix $\ep>0$. Choose $\delta=\delta(\ep)$ such that $m_{d+1}(\partial E_{\delta})<\ep$. Now, choose $\D^*=\D^*(\delta)$ so small that $D_{\n C_s}\subset\partial E_{\delta}$ whenever $\D\leq\D^*$, which can be done since all cells in $\n C_s$ must intersect $\partial E$, and the distance from the furthest point in each such cell to $\partial E$ is at most $\sqrt{dh^2+\tau^2}$. Therefore we need only pick $\D$ so small that $\sqrt{dh^2+\tau^2}<\delta$ to guarantee $D_{\n C_s}\subset\partial E_{\delta}$. It follows that
\[
m_{d+1}(D_{\n C_s})\leq m_{d+1}(\partial E_{\delta})<\ep
\]
for each $\D\leq\D^*$. Therefore $m_{d+1}(D_{\n C_s})\ra0$ as $\D\ra0$. Since $s\in L_2(D)$, it follows by the absolute continuity of the integral that
\[
\Vert s\Vert_{L_2(D_{\n C_s})}\longrightarrow0~~\text{as}~~\D\ra0,
\]
hence from (\ref{calc1}) we conclude $|I_{\D}|\ra0$ as $\D\ra0$.\\

Next, observe that
\begin{gather}
\int\limits_0^T\int\limits_{\Omega}\Big(\tilde V_{\D}^i(x,t)-\tilde V_{\D}^i(x+hz,t)\Big)s(x,t)\,dx\,dt=\nonumber\\[4mm]=I_1+I_2+I_3,\label{calc2}
\end{gather}
where
\[
I_1=\int\limits_0^T\int\limits_{\Omega}\Big(\tilde V_{\D}^i(x,t)-g(x,t)\Big)s(x,t)\,dx\,dt,
\]
\[
I_2=\int\limits_0^T\int\limits_{\Omega}\Big(g(x,t)-g(x+hz,t)\Big)s(x,t)\,dx\,dt,
\]
\[
I_3=\int\limits_0^T\int\limits_{\Omega}\Big(g(x+hz,t)-\tilde V_{\D}^i(x+hz,t)\Big)s(x,t)\,dx\,dt.
\]
We claim each of $|I_1|,|I_2|,|I_3|$ go to $0$ as $\D'\ra0$. Since $g$ is the weak limit of $\tilde V^i_{\D'}$, it follows $|I_1|\ra0$ as $\D'\ra0$. Since $g\in L_2(D)$, by Cauchy-Schwartz inequality and $L_2$-norm continuity of the translation it follows $|I_2|\ra0$ as $h\ra0$. As for $I_3$, by the change of variable $u=x+hz$, we note
\begin{gather}
I_3=\int\limits_0^T\int\limits_{\Omega+hz}\Big(g(u,t)-\tilde V^i_{\D}(u,t)\Big)s(u-hz,t)\,du\,dt=\nonumber\\[4mm]=\int\limits_0^T\int\limits_{\Omega+hz}\Big(g(u,t)-\tilde V^i_{\D}(u,t)\Big)\Big[s(u-hz,t)-s(u,t)+s(u,t)\Big]\,du\,dt=\nonumber\\[4mm]=\int\limits_0^T\int\limits_{\Omega}\Big(g(u,t)-\tilde V^i_{\D}(u,t)\Big)s(u,t)\,du\,dt+\int\limits_0^T\int\limits_{(\Omega+hz)\backslash\Omega}\Big(g(u,t)-\tilde V^i_{\D}(u,t)\Big)s(u,t)\,du\,dt-\nonumber\\[4mm]-\int\limits_0^T\int\limits_{(\Omega\backslash\Omega+hz)}\Big(g(u,t)-\tilde V^i_{\D}(u,t)\Big)s(u,t)\,du\,dt+\nonumber\\[4mm]+\int\limits_0^T\int\limits_{\Omega+hz}\Big(g(u,t)-\tilde V^i_{\D}(u,t)\Big)\Big[s(u-hz,t)-s(u,t)\Big]\,du\,dt=\nonumber\\[4mm]=I_{31}+I_{32}+I_{33}+I_{34}.\label{calc3}
\end{gather}
Through Cauchy-Schwartz, the uniform boundedness of $g-\tilde V^i_{\D}$ in $L_2(\Omega+B_1(0))$, and due to $L_2$-norm continuity of the translation, it follows that $|I_{34}|\ra0$ as $\D\ra0$. Also, $|I_{31}|\ra0$ as $\D'\ra0$ since $\tilde V^i_{\D'}$ converges weakly to $g$ on $D$. $I_{32}\equiv0$ since $s\equiv0$ on $D+hz\backslash D$. $|I_{33}|$ is estimated as follows: apply the Cauchy-Schwartz Inequality, then we have
\begin{align*}
\left|\int\limits_0^T\int\limits_{\Omega\backslash(\Omega+hz)}\Big(g(u,t)-\tilde V^i_{\D}(u,t)\Big)s(u,t)\,du\,dt\right|&\leq\Vert g-\tilde V^i_{\D}\Vert_{L_2(D)}\Vert s\Vert_{L_2(D\backslash(D+hz))}\leq\\[4mm]&\leq C\Vert s\Vert_{L_2(D\backslash(D+hz))}
\end{align*}
where $C$ is a constant independent of $\D$ since $\{\tilde V^i_{\D}\}$ is uniformly bounded in $L_2(D)$. We note that
\[
m_{d+1}(D\backslash(D+hz))\leq Ch^d\longrightarrow0~~\text{as}~~h\ra0
\]
from which, by the absolute continuity of the integral and $s\in L_2(D)$, it follows that $|I_{33}|$ vanishes as $h\ra0$. Hence $|I_3|\ra0$ as $\D'\ra0$.\\

Therefore, for each $z\in\n Z$, (\ref{calc}) and (\ref{calc2}) imply
\[
\left|\sum\limits_{(\gamma,k)\in\n A(\n C_{\D}^D)}\Big(v_{\gamma x_i}(k)-v_{\gamma^z x_i}(k)\Big)\sum\limits_{j=1}^ma_j\int\limits_{C_{\D}^{(\gamma,k)}\cap E_j}W_z(x)\,dy\right|\longrightarrow0~~\text{as}~~\D'\ra0
\]
uniformly with respect to $z\in\n Z$. Using this result, we conclude from (\ref{greater}) that
\[
\left|\int\limits_0^T\int\limits_{\Omega}\left(\tilde V_{\D}^i(x,t)-\frac{\partial}{\partial x_i}V_{\D}(x,t)\right)s(x,t)\,dx\,dt\right|\longrightarrow0~~\text{as}~~\D'\ra0,
\]
which proves that $g=g^*$ in $L_2(D)$ due to the arbitrariness of $s\in\m S_D$. But since $g,g^*$ were arbitrary weak limit points of $\{\tilde V_{\D}^i\},\{\partial V_{\D}/\partial x_i\}$ respectively, it follows $0$ is the unique weak limit of the sequence $\{\tilde V_{\D}^i-\partial V_{\D}/\partial x_i\}$. Statement (f) follows after running the previous argument through all $i=1,2,\ldots,d$.\hfill{$\square$}

\section{Approximation Theorem}\label{approxtheorem}

\begin{theorem}\label{approx} Let $\{[f]_{\D}\}$ be a sequence of discrete control vectors such that there exists $R>0$ for which $[f]_{\D}\in\n F_{\D}^R$ for each $\D$, and such that the sequence of interpolations $\{\n P_{\D}([f]_{\D})\}$ converges weakly to $f$ in $L_2(D)$. Then the sequence of interpolations $\{V_{\D}'\}$ of associated discrete state vectors converges weakly in $W_2^{1,1}(D)$ to $v=v(x,t;f)\in\overset{\circ}{W}{}_2^{1,1}(D)\cap L_{\infty}(D)$, with $v$ the unique weak solution to the Stefan Problem in the sense of Definition \ref{weaksoldef}.
\end{theorem}

\emph{Proof.} From (a) and (b) of Theorem \ref{interp}, it follows that $\{V_{\D}'\}$ is uniformly bounded in $W_2^{1,1}(D)\cap L_{\infty}(D)$. Consequently, $\{V_{\D}'\}$ has a weak limit point in $W_2^{1,1}(D)$. So let $v\in W_2^{1,1}(D)$ be any weak limit point of $\{V_{\D}'\}$ in $W_2^{1,1}(D)$. By the Rellich-Kondrachev Theorem \cite{Nikolski}, it is known that a subsequence of $\{V_{\D}'\}$ converges strongly to $v$ in $L_2(D)$. This allows one to choose a further subsequence of $\{V_{\D}'\}$ which converges pointwise a.e. to $v$ on $D$. Since $\{V_{\D}'\}$ is uniformly bounded in $L_{\infty}(D)$, we have that $v\in L_{\infty}(D)$. Moreover, by construction, $V'_{\D}\equiv0$ on $S$ for each $\D$. Due to $v$ being a weak limit point of $\{V_{\D}'\}$ in $W_2^{1,1}(D)$, it follows that
\[
0=\lim\limits_{\D'\ra0}\Vert v|_S-V_{\D'}'|_S\Vert_{L_2(S)}=\Vert v|_S\Vert_{L_2(S)}
\]
from which we conclude $v|_S=0$. Thus $v\in\overset{\circ}{W}{}_2^{1,1}(D)\cap L_{\infty}(D)$. Henceforth we proceed to show that $v$ satisfies the integral identity (\ref{weaksol}).\\

For simplicity of notation we write the subsequence of $\{V_{\D}'\}$ that converges weakly to $v$ in $W_2^{1,1}(D)$ and pointwise a.e. on $D$ as the whole sequence $\D$. Let $\psi\in\overset{\bullet}{\m C}{}^1(D)$, where $\overset{\bullet}{\m C}{}^1(D)$ be a space of all continuously differentiable functions on $\overline D$ whose support is a positive distance away from $S$ (the lateral boundary of $D$) and from $\Omega\times\{t=T\}$ (the top of the cylinder $D$). Since $D_{\D}\nearrow D$, it follows that there exists $\D^*$ small enough so that $\overline{\text{supp }\psi}\subset D_{\D}$ for all $\D\leq\D^*$. For each $\D\leq\D^*$, define the collection $[\psi]_{\D}=(\psi_{\gamma}^k)$ indexed by $\n A(D_{\D})$ as
\[
\psi_{\gamma}^k:=\psi(x_{\gamma},t_k).
\]
Per our previous remarks, it is clear that for fixed $k$, the collection $\{\psi_{\gamma}^k\}$ is an admissible test collection for the summation identity (\ref{dsveq}). Moreover we remark that independently of the value of $\tau$, we have $\psi_{\gamma}^n=0$ for all $\gamma\in\n A(\Omega_{\D})$. So fix $k=1,\ldots,n$. Let $\eta_{\gamma}:=\tau\psi_{\gamma}^k$ in (\ref{dsveq}). This gives
\begin{equation}\label{dsvplug}
\tau\sum\limits_{\n A}h^d\left[\big(b_n(v_{\gamma}(k))\big)_{\bar t}\psi_{\gamma}^k+\sum\limits_{i=1}^dv_{\gamma x_i}(k)\psi_{\gamma x_i}^k-f_{(\gamma,k)}^{\D}\psi_{\gamma}^k\right]=0
\end{equation}
for each $k=1,2,\ldots,n$. Add up all identities (\ref{dsvplug}) over $k=1,\ldots,n$ to obtain
\begin{equation}\label{dsvsum}
\sum\limits_{k=1}^n\tau\sum\limits_{\n A}h^d\left[\big(b_n(v_{\gamma}(k))\big)_{\bar t}\psi_{\gamma}^k+\sum\limits_{i=1}^dv_{\gamma x_i}(k)\psi_{\gamma x_i}^k-f_{(\gamma,k)}^{\D}\psi_{\gamma}^k\right]=0.
\end{equation}
By summation by parts we observe
\begin{align}
\sum\limits_{k=1}^n\tau\sum\limits_{\n A}h^d\big(b_n(v_{\gamma}(k))\big)_{\bar t}\psi_{\gamma}^k&=\sum\limits_{k=1}^n\sum\limits_{\n A}h^db_n(v_{\gamma}(k))\psi_{\gamma}^k-\sum\limits_{k=1}^n\sum\limits_{\n A}h^db_n(v_{\gamma}(k-1))\psi_{\gamma}^k=\nonumber\\[4mm]&=\sum\limits_{k=1}^n\sum\limits_{\n A}h^db_n(v_{\gamma}(k))\psi_{\gamma}^k -\sum\limits_{k=0}^{n-1}\sum\limits_{\n A}h^db_n(v_{\gamma}(k))\psi_{\gamma}^{k+1}=\nonumber\\[4mm]&=-\sum\limits_{k=1}^{n-1}\tau\sum\limits_{\n A}h^db_n(v_{\gamma}(k))\psi_{\gamma t}^k-\sum\limits_{\n A}h^db_n(\Phi_{\gamma})\psi_{\gamma}^1\label{sumparts},
\end{align}
where $\psi_{\gamma t}^k$ is the forward time difference. Using (\ref{sumparts}) in (\ref{dsvsum}) we can write
\begin{gather}\label{dsvsum1}
-\sum\limits_{k=1}^{n-1}\tau\sum\limits_{\n A}h^db_n(v_{\gamma}(k))\psi_{\gamma t}^k+\sum\limits_{k=1}^{n}\tau\sum\limits_{\n A}h^d\left[\sum\limits_{i=1}^dv_{\gamma x_i}(k)\psi_{\gamma x_i}^k-f_{(\gamma,k)}^{\D}\psi_{\gamma}^k\right]-\nonumber\\[4mm]-\sum\limits_{\n A}h^db_n(\Phi_{\gamma})\psi_{\gamma}^1=0.
\end{gather}
Define the following interpolations $\overline\Phi_{\D},\overline\psi_{\D},\overline\psi_{\D}^t,\overline\psi_{\D}^i,~i=1,\ldots,d$ of the collections $\{\Phi_{\gamma}\}, [\psi]_{\D}$ and the forward differences of the latter:
\begin{gather}
\overline\Phi_{\D}\Big|_{R_{\D}^{\gamma}}=\Phi_{\gamma},~~\alpha\in\n A,\qquad\overline\Phi_{\D}\equiv0\text{ elsewhere on } D,\nonumber\\[4mm]
\overline\psi_{\D}\Big|_{C_{\D}^{\alpha}}=\psi_{\gamma}^k,~~\alpha\in\n A(\n C_{\D}^D),\qquad\overline\psi_{\D}\equiv0\text{ elsewhere on }D,\nonumber\\[4mm]
\overline\psi_{\D}^t\Big|_{C_{\D}^{\alpha}}=\psi_{\gamma t}^k,~~\alpha\in\n A(\n C_{\D}^D\backslash\n R_{\D}^{\gamma,n}),\qquad\overline\psi_{\D}^t\equiv0\text{ elsewhere on }D,\nonumber\\[4mm]
\overline\psi_{\D}^i\Big|_{C_{\D}^{\alpha}}=\psi_{\gamma x_i}^k,~~\gamma\in\n A,~~k=1,\ldots,n,\qquad\overline\psi_{\D}^i\equiv0\text{ elsewhere on }D.\nonumber
\end{gather}
With these functions and with the interpolations described in Section~\ref{Interpolations},
 identity (\ref{dsvsum1}) can be written in the following way:
\begin{gather}
-\sum\limits_{k=1}^{n-1}~\int\limits_{t_{k-1}}^{t_k}\sum\limits_{\n A}\int\limits_{R_{\D}^{\gamma}}b_n(\tilde V_{\D})\overline\psi_{\D}^t\,dx\,dt~~+~~\sum\limits_{k=1}^{n}~\int\limits_{t_{k-1}}^{t_k}\sum\limits_{\n A}\int\limits_{R_{\D}^{\gamma}}\sum\limits_{i=1}^d\tilde V_{\D}^i\overline\psi_{\D}^i\,dx\,dt~~-\nonumber\\[4mm]-\sum\limits_{\n A(\n C_{\D}^D)}\tau h^d\left(\frac1{\tau h^d}\int\limits_{C_{\D}^{\alpha}}f^{\D}\,dx\,dt\right)\psi_{\gamma}^k~~-~~\sum\limits_{\n A}\int\limits_{R_{\D}^{\gamma}}b_n(\overline\Phi_{\D})\overline\psi_{\D}(x,\tau)\,dx=0,
\nonumber\\[4mm]
\intertext{whence it follows}
\int\limits_0^T\int\limits_{\Omega}\Big[-b_n(\tilde V_{\D})\overline\psi_{\D}^t+\sum\limits_{i=1}^d\tilde V_{\D}^i\overline\psi_{\D}^i-f^{\D}\overline\psi_{\D}\Big]\,dx\,dt-\int\limits_{\Omega}b_n(\overline\Phi_{\D})\overline\psi_{\D}(x,\tau)\,dx=0,\label{int1}
\end{gather}
since $\overline\psi_{\D}^t\equiv0$ on $\Omega\times(T-\tau,T]$.\\

Next we show that the sequences $\{b_n(\tilde V_{\D})\},\{b_n(\overline\Phi_{\D})\}$ converge weakly in $L_2(D)$ to functions of type $\n B$. Due to Theorem \ref{interp} (c),(d), we know that $\{\tilde V_{\D}\}$ converges strongly to $v$ in $L_2(D)$. As such, we can extract a subsequence of $\{\tilde V_{\D}\}$ that converges pointwise a.e. on $D$ to $v$. For ease of notation let this subsequence be denoted as the whole sequence. Define the set
\[
N:=\Big\{(x,t)\in D~\Big|~\lim\limits_{\D\ra0}|\tilde V_{\D}(x,t)-v(x,t)|\neq0.\Big\},
\]
and from the previous remarks it's clear $m_{d+1}(N)=0$. Now fix arbitrary $(x,t)\in D\backslash N$. For such $(x,t)$, we have
\[
\tilde V_{\D}(x,t)\longrightarrow v(x,t),~~\text{as}~~\D\ra0.
\]
Suppose that at the point $(x,t)\in D\backslash N$ we have $v(x,t)\neq v^j$ for any $j=1,\ldots,J$ (recall the $v^j$'s correspond to phase transition temperatures). In this case we observe
\[
b_n(\tilde V_{\D}(x,t)) = \int\limits_{\tilde V_{\D}(x,t)-\frac1n}^{\tilde V_{\D}(x,t)+\frac1n}\omega_{1/n}(|\tilde V_{\D}(x,t)-u|)b(u)\,du~\longrightarrow~b(v(x,t)),~~\text{as}~~\D\ra0.
\]
On the contrary, if at the point $(x,t)\in D\backslash N$ we have $v(x,t)=v^j$ for some $j\in\{1,\ldots,J\}$, then we have
\[
b(v^j)^- \leq \liminf_{\D\rightarrow \infty}b_n(\tilde V_{\D}(x,t)) \leq \limsup_{\D\rightarrow \infty}b_n(\tilde V_{\D}(x,t)) \leq b(v^j)^+.
\]
The past few observations show that we can pass to a subsequence of $\{b_n(\tilde V_{\D})\}$ which converges pointwise on $D\backslash N$ to a function $\tilde b(x,t)$ that satisfies
\begin{align*}
&\tilde b(x,t)=b(v(x,t))~~\text{whenever}~~v(x,t)\neq v^j,~~\text{and}\\&\tilde b(x,t)\in[b(v^j)^-,b(v^j)^+]~~\text{whenever}~~v(x,t)=v^j\text{ for some }j,
\end{align*}
which shows that $\tilde b$ is a function of type $\n B$ as in Definition \ref{typeB}. Moreover we claim that $\{b_n(\tilde V_{\D})\}$ converges weakly in $L_2(D)$ to $\tilde b$. To see this, it is enough to show that $\tilde b\in L_2(D)$ and that $\{b_n(\tilde V_{\D})\}$ is uniformly bounded in $L_2(D)$. Let $\n V_{\D}$ be the range of $\tilde V_{\D}$. Due to Theorem \ref{boundedness}, it follows the set $\n V=\cup_{\D}\n V_{\D}$ is bounded in $\bb R$, hence its closure $\overline{\n V}$ is compact in $\bb R$. Because of the piecewise continuity of $b$, the sequence $\{b(\tilde V_{\D}(x,t))\}$ is uniformly bounded in $L_{\infty}(D)$, and so too must be the sequence $\{b_n(\tilde V_{\D})\}$. Hence $\{b_n(\tilde V_{\D})\}$ is uniformly bounded in $L_2(D)$ as well, since $D$ is a set of finite measure. A very similar argument concludes that $\tilde b\in L_2(D)$ too.

We have proved that a subsequence of $\{b_n(\tilde V_{\D})\}$ converges weakly in $L_2(D)$ to $\tilde b$, a function of type $\n B$. It is proved in a completely analogous way that a further subsequence of $\{b_n(\overline\Phi_{\D})\}$ converges weakly to $\tilde b_0$, a function of type $\n B$. Again we denote this further subsequence as the whole sequence, for simplicity of notation.\\

Carrying on, it is easily shown that the functions $\overline\psi_{\D},\overline\psi_{\D}^t,\overline\psi_{\D}^i$ converge uniformly on $\overline D$ to the functions $\psi,\partial\psi/\partial t,\partial\psi/\partial x_i$ respectively as $\D\ra0$. Consequently, (\ref{int1}) implies
\begin{equation}\label{int2}
\int\limits_0^T\int\limits_{\Omega}\Big[-b_n(\tilde V_{\D})\frac{\partial\psi}{\partial t}+\sum\limits_{i=1}^d\tilde V_{\D}^i\frac{\partial\psi}{\partial x_i} -f^{\D}\psi\Big]\,dx\,dt-\int\limits_{\Omega}b_n(\overline\Phi_{\D})\psi(x,0)\,dx+I=0,
\end{equation}
where
\begin{gather}
I=\int\limits_0^T\int\limits_{\Omega}\Big[-b_n(\tilde V_{\D})\left(\overline\psi_{\D}^t-\frac{\partial\psi}{\partial t}\right)+\sum\limits_{i=1}^d\tilde V_{\D}^i\left(\overline\psi_{\D}^i- \frac{\partial\psi}{\partial x_i}\right) -f^{\D}\big(\overline\psi_{\D}-\psi\big)\Big]\,dx\,dt-\nonumber\\[4mm]-\int\limits_{\Omega}b_n(\overline\Phi_{\D})\Big(\overline\psi_{\D}(x,\tau)-\psi(x,0)\Big)\,dx.\label{I}
\end{gather}
We claim $|I|\ra0$ as $\D\ra0$. Since the sequences $\{b_n(\tilde V_{\D})\},\{\tilde V_{\D}^i\},\{f^{\D}\}$ are uniformly bounded in $L_2(D)$, and since $\overline\psi_{\D},\overline\psi_{\D}^t,\overline\psi_{\D}^i$ converge uniformly on $\overline D$ to the functions $\psi,\partial\psi/\partial t,\partial\psi/\partial x_i$ respectively as $\D\ra0$ (hence, strongly in $L_2(D)$), then by the Cauchy-Schwartz inequality it is seen that the absolute value of the $D-$integral term of (\ref{I}) vanishes as $\D\ra0$. As for the last term, we observe
\begin{align*}
\left|\int\limits_{\Omega}b_n(\overline\Phi_{\D})\Big(\overline\psi_{\D}(x,\tau)-\psi(x,0)\Big)\,dx\right|&\leq\Vert b_n(\overline\Phi_{\D})\Vert_{L_2(\Omega)}~\Vert\overline\psi_{\D}(x,\tau)-\psi(x,0)\Vert_{L_2(\Omega)}\leq\\[4mm]&\leq C\Big(\Vert\overline\psi_{\D}(x,\tau)-\psi(x,\tau)\Vert_{L_2(\Omega)}+\Vert\psi(x,\tau)-\psi(x,0)\Vert_{L_2(\Omega)}\Big),
\end{align*}
and both terms on the right-hand side of the above inequality converge to $0$ as $\D\ra0$ (the first due to uniform convergence of $\{\overline\psi_{\D}\}$ to $\psi$ on $D$, and the second due to uniform continuity of $\psi$). Therefore $|I|\ra0$ as $\D\ra0$. So, due to the weak convergence of the sequences $\{b_n(\tilde V_{\D})\},\{\tilde V_{\D}^i\},\{f^{\D}\}, \{b_n(\overline\Phi_{\D})\}$ to the functions $\tilde b(x,t),\partial v/\partial x_i, f, \tilde b_0$ in $L_2(D)$ and $L_2(\Omega)$ respectively, it follows that taking $\D\ra0$ on (\ref{int2}) gives the identity
\begin{equation}
\int\limits_0^T\int\limits_{\Omega}\Big[-\tilde b(x,t)\frac{\partial\psi}{\partial t}+\sum\limits_{i=1}^d\frac{\partial v}{\partial x_i}\frac{\partial\psi}{\partial x_i} -f\psi\Big]\,dx\,dt-\int\limits_{\Omega}\tilde b_0(x)\psi(x,0)\,dx=0.\nonumber
\end{equation}
which is (\ref{weaksol}). Thus we have proved $v$ satisfies integral identity (\ref{weaksol}) for some functions $b,b_0$ of type $\n B$, and for arbitrary test function $\psi\in\overset{\bullet}{\m C}{}^1(D)$. Since  $\overset{\bullet}{\m C}{}^1(D)$ is dense in the set of admissible test functions for integral identity (\ref{weaksol}) and due to Remark \ref{some}, we have that $v$ is a weak solution to the Stefan Problem in the sense of Definition \ref{weaksoldef}.
Therefore, we have proved that if $v$ is a weak limit point of $\{V_{\D}'\}$, then it must be a weak solution to the Stefan Problem. Due to uniqueness of the weak solution \cite{LSU}
(see Remark~\ref{some}) it follows that $\{V_{\D}'\}$ has one and only one weak limit point, which shows that the whole sequence $\{V_{\D}'\}$ converges weakly to $v$ in $W_2^{1,1}(D)$. This ends the proof of the theorem.\hfill{$\square$}\\

Theorem \ref{approx} readily provides us with a general existence theorem for the Stefan Problem.
\begin{corollary}\label{vestimates} Let $f\in L_{\infty}(D)$. Then there exists $v=v(x,t;f)\in\overset{\circ}{W}{}_2^{1,1}(D)\cap L_{\infty}(D)$ which is a weak solution to the Stefan Problem. Moreover, $v$ satisfies the following estimates:
\begin{equation}\label{boundv}
\Vert v\Vert_{L_{\infty}(D)}\leq e^T\max\left\{\frac1{\bar b}\Vert f\Vert_{L_{\infty}(D)}~,~\Vert\Phi\Vert_{L_{\infty}(\Omega)}\right\},
\end{equation}
\begin{equation}\label{energyv}
\Vert D_xv\Vert_{L_2(D)}^2+\Vert v_t\Vert_{L_2(D)}^2\leq C\left[~\Vert f\Vert_{L_2(D)}^2+\Vert\Phi\Vert_{W_2^1(\Omega)}^2\right]
\end{equation}
where $C$ is a constant depending on $\bar b$ and $d$.
\end{corollary}

\emph{Proof.} Given $f\in L_{\infty}(D)$, consider the collection $[f]_{\D}:=\n Q_{\D}(f)$. Then the interpolations $\n P_{\D}([f]_{\D})$ converge strongly to $f$ in $L_2(D)$, and by Lemma \ref{PQ} and Cauchy-Schwartz inequality we have
\[  \Vert[f]_{\D}\Vert_{\ell_{\infty}}\leq\Vert f\Vert_{L_{\infty}(D)}, \ \ \Vert f^{\D}\Vert_{L_2(D)}\leq \Vert f\Vert_{L_2(D)}. \]
The conditions of Theorem \ref{approx} are satisfied, so there exists $v=v(x,t;f)\in\overset{\circ}{W}{}_2^{1,1}(D)\cap L_{\infty}(D)$ which is a weak solution to the Stefan Problem in the sense of Definition \ref{weaksoldef}. Moreover, the sequence $\{V_{\D}'\}$ converges to $v$ weakly in $W_2^{1,1}(D)$, and strongly in $L_2(D)$. In particular, there is a subsequence which converges to $v$ almost everywhere on $D$. By Theorem \ref{boundedness} and Theorem \ref{interp}(a), it is clear that for each $\D$, $\Vert V_{\D}'\Vert_{L_{\infty}(D)}$ is bounded above by the right-hand side of (\ref{boundedest}). Therefore, (\ref{boundv}) easily follows. Furthermore, we have
\begin{equation}\label{cor2}
\Vert D_xv\Vert_{L_2(D)}\leq\liminf_{\D\ra0}\Vert D_xV_{\D}'\Vert_{L_2(D)}, \ \  \Vert v_t\Vert_{L_2(D)}\leq\liminf_{\D\ra0}\Big\Vert \frac{\partial}{\partial t}V_{\D}'\Big\Vert_{L_2(D)}.
\end{equation}
Using the estimations (\ref{pwlineartok}), (\ref{pwlinearxok}), from (\ref{cor2}), (\ref{energyv}) follows.\hfill{$\square$}

\section{Existence of the Optimal Control}\label{exists}

\emph{Proof of Theorem \ref{optsol}.} By definition of $\n J_*$, there exists a sequence $\{f_{\ell}\}\subset\n F^R$ such that $\n J(f_{\ell})\searrow\n J_*$. Such a sequence is uniformly bounded in $L_2(D)$ since $D$ is bounded, so the sequence has a weak limit point $f$ in $L_2(D)$. We claim $f\in\n F^R$. By Mazur's Lemma, there is a sequence $\{F_{\ell}\}$ given as
\[
F_{\ell}(x,t)=\sum\limits_{k=\ell}^{K(\ell)}a_k^{\ell}f_k(x,t)
\]
which converges strongly to $f$ in $L_2(D)$ as $\ell\ra\infty$, where for each $\ell$, the set of real numbers $\{a_{\ell}^{\ell},\ldots,a_{K(\ell)}^{\ell}\}$ is contained in $[0,1]$ and
\[
\sum\limits_{k=\ell}^{K(\ell)}a_k^{\ell}=1.
\]
Then there is a subsequence $F_{\ell_m}$ which converges pointwise a.e. on $D$ to $f$ as $m\ra\infty$. We observe that
\[
\Vert F_{\ell}\Vert_{L_{\infty}(D)}\leq\sum\limits_{k=\ell}^{K(\ell)}a_k^{\ell}\Vert f_k\Vert_{L_{\infty}(D)}\leq R
\]
uniformly over $\ell$. Therefore, it follows that $f\in\n F^R$.\\

Corollary \ref{vestimates} implies the existence of the unique weak solutions to the Stefan Problem for any of the functions $f_{\ell},f$. So let $v_{\ell}=v(x,t;f_{\ell}), v=v(x,t;f)$ be the unique weak solutions to the Stefan Problem with $f_{\ell}$ and $f$ as controls, respectively. Due to (\ref{boundv}), (\ref{energyv}) and the fact that $\Vert f_{\ell}\Vert_{L_{\infty}(D)}\leq R$ for all $\ell$, it follows that the sequence $\{v_{\ell}\}$ is uniformly bounded in the spaces $W_2^{1,1}(D)$ and $L_{\infty}(D)$. Therefore, $\{v_{\ell}\}$ has a weak limit point in $W_2^{1,1}(D)$. Let $\tilde v$ be such a weak limit point, and for ease of notation say that the whole sequence $\{v_{\ell}\}$ converges to $\tilde v$ weakly in $W_2^{1,1}(D)$. It's clear then that $\tilde v\in W_2^{1,1}(D)\cap L_{\infty}(D)$. Moreover, since $v_{\ell}|_S=0$
for each $\ell$, it follows that $\tilde v|_S=0$. Hence $\tilde v\in\overset{\circ}{W}{}_2^{1,1}(D)\cap L_{\infty}(D)$.\\

Next we show that $\tilde v$ is actually a weak solution to the Stefan Problem with $f$ as control. To this end, fix an arbitrary $B(x,t,v)$, a function of type $\n B$, and $\psi$ an arbitrary admissible test function for integral identity (\ref{weaksol}). For each $\ell\in\bb N$, we have the identity
\begin{equation}\label{ellv}
\int\limits_D\Big[-B(x,t,v_{\ell}(x,t))\psi_t+\nabla v_{\ell}\cdot\nabla\psi-f_{\ell}\psi\Big]\,dxdt - \int\limits_{\Omega}B_0(x,0,\Phi(x))\psi(x,0)\,dx=0.
\end{equation}
Since $\{v_{\ell}\},\{f_{\ell}\}$ converge weakly to $\tilde v,f$ in $W_2^{1,1}(D)$ respectively, to obtain the desired identity it is only left to show that
\begin{equation}\label{Bweak}
B(x,t;v_{\ell}(x,t))\longrightarrow B'(x,t;\tilde v(x,t))~~\text{weakly in }L_2(D)~~\text{as }\ell\ra\infty
\end{equation}
where $B'$ is some function of type $\n B$. To see that (\ref{Bweak}) is true, first pass to a subsequence $\{v_{\ell_m}\}$ that converges pointwise a.e. on $D$ to $\tilde v$, and for ease of notation write this subsequence as the whole sequence. It is sufficient to prove that that there exists some function $B'$ of type $\n B$ such that
\begin{itemize}
\item[] (i) $B(x,t;v_{\ell}(x,t))\longrightarrow B'(x,t;\tilde v(x,t))~~\text{pointwise a.e. on }D~~\text{as }\ell\ra\infty$,
\item[] (ii)  $B(x,t;v_{\ell}(x,t))$ is uniformly bounded in $L_2(D)$, 
\item[] (iii) $B'(x,t;\tilde v(x,t))\in L_2(D)$ 
\end{itemize}
To prove (i), let
\[
N:=\Big\{(x,t)\in D~\Big|~\lim\limits_{\ell\ra0}|v_{\ell}(x,t)-\tilde v(x,t)|\neq0.\Big\}.
\]
Then by construction, $m_{d+1}(N)=0$ and $\{v_{\ell}\}$ converges pointwise to $\tilde v$ on $D\backslash N$. Now fix $(x,t)\in D\backslash N$. Suppose that
\[
\tilde v(x,t)\neq v^j~~\text{for any }j.
\] 
In this case we note that $b$ is continuous at $\tilde v(x,t)$, and therefore
\[
B(x,t;v_{\ell}(x,t))=b(v_{\ell}(x,t))\longrightarrow b(\tilde v(x,t))=B(x,t;\tilde v(x,t))
\]
as $\ell\ra\infty$. On the contrary, suppose that
\[
\tilde v(x,t)=v^j~~\text{for some }j\in\{1,2,\ldots,J\}.
\]
By way of contradiction, assume that there is a subsequence $\{B(x,t;v_{\ell_m}(x,t))\}$ such that
\[
L:=\lim\limits_{m\ra\infty}B(x,t;v_{\ell_m}(x,t))\notin[b(v^j)^-,b(v^j)^+].
\]
Then dist$\big(L,[b(v^j)^-,b(v^j)^+]\big)>0$. Since $b$ is monotone, this gives a contradiction to the fact that $v_{\ell_m}(x,t)\ra v^j$. Thus in this case we have
\[
b(v^j)^-\leq\liminf\limits_{\ell\ra\infty}B(x,t;v_{\ell}(x,t))\leq\limsup\limits_{\ell\ra\infty}B(x,t;v_{\ell}(x,t))\leq b(v^j)^+.
\]
Hence the assertion (i) is proved. Since $\tilde v\in L_{\infty}(D)$ and $\{v_{\ell}\}$ is uniformly bounded in $L_{\infty}(D)$, the assertions (ii) and (iii) easily follow from the definition of the functional class $\n B$. Therefore (\ref{Bweak}) is true, and so passing $\ell\ra\infty$ on (\ref{ellv}), we obtain
\[
\int\limits_D\Big[-B'(x,t,\tilde v(x,t))\psi_t+\nabla \tilde v\cdot\nabla\psi-f\psi\Big]\,dxdt - \int\limits_{\Omega}B_0(x,0,\Phi(x))\psi(x,0)\,dx=0,
\]
from which we conclude that $\tilde v$ is a weak solution to the Stefan Problem with $f$ as a control. Due to uniqueness, we then have $\tilde v$ is the same element as $v$ in the space $\overset{\circ}{W}{}_2^{1,1}(D)\cap L_{\infty}(D)$.\\

By employing the following elementary identity for elements $a,b,c$ of the Hilbert space $H$
\begin{gather}
\Vert a-b\Vert_H^2-\Vert c-b\Vert_H^2=
\langle a-c,a-c\rangle-2\langle a-c,b-c\rangle,\label{Hcomp}
\end{gather}
we have
\begin{gather}
|\n J(f_{\ell})-\n J(f)|=\Big|\Vert v_{\ell}|_{\Omega\times\{t=T\}}-\Gamma\Vert^2_{L_2(\Omega)}-\Vert v|_{\Omega\times\{t=T\}}-\Gamma\Vert^2_{L_2(\Omega)}\Big|=\nonumber\\[4mm]=\left|\Vert v_{\ell}|_{\Omega\times\{t=T\}}-v|_{\Omega\times\{t=T\}} \Vert_{L_2(\Omega)}^2-2\int\limits_{\Omega}\Big(v_{\ell}|_{\Omega\times\{t=T\}}-v|_{\Omega\times\{t=T\}}\Big)(\Gamma(x)-v|_{\Omega\times\{t=T\}})\,dx\right|.\label{Jdiff}
\end{gather}
Since we've shown that $v_{\ell}\ra v$ weakly in $W_2^{1,1}(D)$, and since weak convergence in $W_2^{1,1}(D)$ implies strong convergence in the space of traces, it follows that $v_{\ell}|_{\Omega\times\{t=T\}}\ra v|_{\Omega\times\{t=T\}}$ strongly in $L_2(\Omega)$. As a result, (\ref{Jdiff}) implies that
\[
\n J(f_{\ell})\ra\n J(f)~~\text{as}~~\ell\ra\infty,
\]
so that $\n J(f)=\n J_*$. Theorem is proved.\hfill{$\square$}

In the previous proof we have actually shown the
\begin{corollary}\label{jweakcont} The cost functional $\n J$ is weakly continuous on $\n F^R$ for any $R>0$.
\end{corollary}

\section{Convergence of the Discrete Optimal Control Problem}\label{approxfunc}

\emph{Proof of Theorem \ref{funcapprox}.} To prove (\ref{approx1}) and (\ref{approx2}), it is enough to show that conditions (i) and (ii) of Lemma \ref{Vasil} are satisfied. We first claim that for any $f\in\n F^R$,
\begin{equation}\label{itoj}
\lim\limits_{\D\ra0}\left|\n I_{\D}(\n Q_{\D}(f))-\n J(f)\right|=0.
\end{equation}
To see this, first note that if we write $[f]_{\D}=\n Q_{\D}(f)$, then
\begin{equation}\label{convweak}
\n P_{\D}([f]_{\D})\longrightarrow f~~\text{weakly in }L_2(D)\text{ as }\D\ra0,
\end{equation}
whence we have by Theorem \ref{approx} that the interpolations $V_{\D}'$ of the discrete state vectors $[v([f]_{\D})]_{\D}$ converge weakly in $W_2^{1,1}(D)$ to the unique weak solution $v=v(x,t;f)$ of the Stefan problem with control $f$ (the convergence in (\ref{convweak}) can be taken to be strong, but the argument given here only assumes weak convergence). Consequently, it is known that this implies
\begin{equation}\label{traceconv}
V_{\D}'\Big|_{\Omega\times\{t=T\}}\longrightarrow v\Big|_{\Omega\times\{t=T\}}~~\text{strongly in }L_2(\Omega\times\{t=T\})~~\text{as}~~\D\ra0.
\end{equation}

Define $\tilde\Gamma_{\D}$ as the piece-wise constant interpolation of the collection $\{\Gamma_{\gamma}\}$:
\[
\tilde\Gamma_{\D}\Big|_{R_{\D}^{\gamma}}=\Gamma_{\gamma},~~\forall\gamma\in\n A,\qquad\tilde\Gamma_{\D}\equiv0~~\text{elsewhere on }\Omega.
\]
Next, we note
\begin{align}
\n I_{\D}(\n Q_{\D}(f))=\sum\limits_{\n A}h^d|v_{\gamma}(n)-\Gamma_{\gamma}|^2=\sum\limits_{\n A}\int\limits_{R_{\D}^{\gamma}}\Big|\tilde V_{\D}\Big|_{t=T}-\tilde\Gamma_{\D}\Big|^2\,dx=\left\Vert\tilde V_{\D}\Big|_{t=T}-\tilde\Gamma_{\D}\right\Vert_{L_2(\Omega)}^2\nonumber,
\end{align}
so, using (\ref{Hcomp}), we observe
\begin{gather}
\n I_{\D}(\n Q_{\D}(f))-\n J(f)=\nonumber\\[4mm]=\left\Vert\tilde V_{\D}\Big|_{t=T}-\tilde\Gamma_{\D}\right\Vert_{L_2(\Omega)}^2-\left\Vert v\Big|_{t=T}-\tilde\Gamma_{\D}\right\Vert_{L_2(\Omega)}^2+\left\Vert v\Big|_{t=T}-\tilde\Gamma_{\D}\right\Vert_{L_2(\Omega)}^2-\left\Vert v\Big|_{t=T}-\Gamma\right\Vert^2_{L_2(\Omega)}=\nonumber\\[4mm]=\left\Vert\tilde V_{\D}\Big|_{t=T}-v\Big|_{t=T}\right\Vert_{L_2(\Omega)}^2-2\left\langle\tilde V_{\D}\Big|_{t=T}-v\Big|_{t=T}~,~\tilde\Gamma_{\D}-v\Big|_{t=T}\right\rangle_{L_2(\Omega)}+\nonumber\\[4mm]+\left\Vert\tilde \Gamma_{\D}-\Gamma\right\Vert_{L_2(\Omega)}^2-2\left\langle\tilde \Gamma_{\D}-\Gamma~,~v\Big|_{t=T}-\Gamma\right\rangle_{L_2(\Omega)}.\label{itoj1}
\end{gather}
By convergence of the Steklov averages to the original function in $L_2$, it follows that
\begin{equation}\label{steklovconvergence}
\left\Vert\tilde \Gamma_{\D}-\Gamma\right\Vert_{L_2(\Omega)}\longrightarrow0~~\text{as}~~h\ra0.
\end{equation}
We also have
\begin{align}
\left\Vert\tilde V_{\D}\Big|_{t=T}-v\Big|_{t=T}\right\Vert_{L_2(\Omega)}&\leq\left\Vert\tilde V_{\D}\Big|_{t=T}-V_{\D}'\Big|_{t=T}\right\Vert_{L_2(\Omega)}+\left\Vert V_{\D}'\Big|_{t=T}-v\Big|_{t=T}\right\Vert_{L_2(\Omega)}\leq\nonumber\\[4mm]&\leq\left\Vert\tilde V_{\D}\Big|_{t=T}-V_{\D}^n\right\Vert_{L_2(\Omega)}+\left\Vert V_{\D}'\Big|_{t=T}-v\Big|_{t=T}\right\Vert_{L_2(\Omega)}\ra 0,\label{2}
\end{align}
as $\D\ra0$. The latter follows from Theorem \ref{interp}(d) and (\ref{traceconv}). By using (\ref{steklovconvergence}), (\ref{2}) and  the uniform boundedness of $\{\tilde\Gamma_{\D}\}$ in $L_2(\Omega)$ from (\ref{itoj1}) it follows that
\begin{equation}\label{willimply}
\left|\n I_{\D}(\n Q_{\D}(f))-\n J(f)\right|\leq C\left(\left\Vert\tilde V_{\D}\Big|_{t=T}-v\Big|_{t=T}\right\Vert_{L_2(\Omega)}+\left\Vert\tilde \Gamma_{\D}-\Gamma\right\Vert_{L_2(\Omega)}\right),
\end{equation}
where $C$ is a constant independent of $\D$. Hence, (\ref{itoj}) is proved.

Next we claim that for any sequence $\{[f]_{\D}\}$ of discrete controls such that $[f]_{\D}\in\n F_{\D}^R$ for some fixed $R>0$, it follows that
\begin{equation}\label{last}
\lim\limits_{\D\ra0}\Big|\n J(\n P_{\D}([f]_{\D}))-\n I_{\D}([f]_{\D})\Big|=0.
\end{equation}
To this end, notice that by Proposition \ref{PQ}, the sequence $\{\n P_{\D}([f]_{\D})\}$ is uniformly bounded in $L_{\infty}(D)$, hence also in $L_2(D)$. Therefore, there is an $L_2(D)$-weak limit point to the sequence $\{\n P_{\D}([f]_{\D})\}$. So let $f\in L_2(D)$ be any weak limit point of the aforementioned sequence, and pass to a subsequence that converges to it in the $L_2(D)-$weak sense. For ease of notation, denote the subsequence as the whole sequence. Then we see that (\ref{convweak}) is true, so the argument leading to the proof of (\ref{itoj}) gives us
\begin{equation}\label{itoj2}
\lim\limits_{\D\ra0}\Big|\n I_{\D}([f]_{\D})-\n J(f)\Big|=0.
\end{equation}
Equipped with (\ref{itoj2}) and Corollary \ref{jweakcont}, we deduce
\[
\Big|\n J(\n P_{\D}([f]_{\D}))-\n I_{\D}([f]_{\D})\Big|\leq\Big|\n J(\n P_{\D}([f]_{\D}))-\n J(f)\Big|+\Big|\n I_{\D}([f]_{\D})-\n J(f)\Big|\longrightarrow0~~\text{as}~~\D\ra0,
\]
from which (\ref{last}) follows, since $f$ was any weak limit point in $L_2(D)$ of $\{\n P_{\D}([f]_{\D})\}$.\\

The results (\ref{itoj}) and (\ref{last}) show that conditions (i) and (ii) of Lemma \ref{Vasil} are satisfied. Therefore (\ref{approx1}) and (\ref{approx2}) follow. Now let $[f]_{\D,\ep}\in\n F_{\D}^R$ be a sequence satisfying (\ref{approxcond}). It is clear that $\{\n P_{\D}([f]_{\D})\}$ is uniformly bounded in $L_2(D)$. Let $f_*$ be any weak limit point of $\{\n P_{\D}([f]_{\D})\}$ in $L_2(D)$. By Corollary \ref{jweakcont} and (\ref{approx2}), we easily see $\n J(f_*)=\n J_*$, hence $f_*\in\n F_*$. The rest of the theorem is an easy consequence of Theorem \ref{approx}.\hfill{$\square$}

\end{document}